\documentclass[12pt]{amsart}
\usepackage[margin=3cm]{geometry}
\usepackage[foot]{amsaddr}
\usepackage{amsmath,amsthm,amssymb}
\usepackage{enumerate}
\usepackage[T1]{fontenc}
\usepackage[utf8]{inputenc}
\usepackage{pgfplots}
\usepackage{subcaption}
\usepackage{float}
\usepackage{dsfont}
\usepackage{comment}
\usepackage{xcolor}

\newtheorem{theorem}{Theorem}[section]
\newtheorem{corollary}[theorem]{Corollary}

\newtheorem{lemma}[theorem]{Lemma}
\newtheorem{proposition}{Proposition}

\theoremstyle{definition}

\newtheorem{remark}{Remark}

\newcommand{\R}{\mathbb R}

\numberwithin{equation}{section}

\begin{document}
\title[Evolution of anisotropic diffusion]{Evolution of anisotropic diffusion in
two-dimensional heterogeneous environments}


\author{Emeric Bouin, Guillaume Legendre, Yuan Lou, Nichole Slover}
\address{CEREMADE, UMR CNRS 7534, Universit\'e Paris-Dauphine, Universit\'e PSL, Place du Mar\'echal De Lattre De Tassigny, 75775 Paris cedex 16, France}
\curraddr{}
\email{bouin@ceremade.dauphine.fr, guillaume.legendre@ceremade.dauphine.fr}

\address{
Department of Mathematics, Ohio State University, Ohio, 43210 USA}
\curraddr{}
\email{lou@math.ohio-state.edu, slover.9@osu.edu}


\date{\today}

\dedicatory{}

\begin{abstract}
We consider a system of two competing populations in two-dimensional heterogeneous environments. The populations are assumed to move horizontally and vertically with different probabilities, 
but are otherwise identical. We regard these probabilities as dispersal strategies.
We show that 
the evolutionarily stable strategies are to move in one direction only. 
Our results predict that it is more beneficial for the species to choose the direction with smaller variation in the resource distribution. This finding seems to be in agreement with the classical results of Hasting \cite{Hastings1983} and Dockery \textit{et al.} \cite{DHMP} for the evolution of slow dispersal, \textit{i.e.} random diffusion is selected against in spatially heterogeneous  environments. These conclusions also suggest that broader dispersal strategies should be considered regarding the movement in heterogeneous habitats.
\end{abstract}

\maketitle

\section{Introduction} 

\subsection {Background and motivation}
In this paper, we consider populations of individuals that disperse in a bounded two-dimensional habitat, where the resources are distributed heterogeneously across the landscape. A natural question is how organisms should distribute themselves in the space to better match the available resources and, accordingly, what kind of dispersal strategies organisms should adopt to reach such distributions; see, \textit{e.g.} \cite{Clobert2001, Fretwell}. Most of previous studies on the evolution of dispersal assume that individuals move in two orthogonal directions with the same probability, which we refer as random dispersal; see, \textit{e.g.} \cite{Hastings1983}. 
For spatially varying but temporally constant environments, Hastings considered a scenario in which the resident is randomly dispersing and reaches the equilibrium; some rare mutant, which is also randomly dispersing but differs from the resident only in the diffusion rate, is introduced to the habitat. He found that slower rates of dispersal will be selected, as the mixing of populations tends to reduce the growth \cite{Altenberg1,Altenberg2}. Dockery \textit{et al.} \cite{DHMP} considered a system of two randomly diffusing competing populations in spatially varying but temporally constant environments, and two species are identical except their dispersal rates.  They showed that the population with the smaller dispersal rate always drive the population with the larger dispersal rate to extinction, irrelevant of the initial data.  This phenomenon is often termed as the evolution of slow dispersal, as any population with a positive dispersal rate will be replaced by a mutant with a smaller dispersal rate.

If we assume that individuals move, \textit{e.g.} horizontally and vertically with two different probabilities (with the sum of probabilities equal to one)  and regard these probabilities as dispersal strategies, what kind of strategies will be evolutionarily stable? Here, we are referring to the situation where the population moves east {or} west with probability $p/2$, and north {or} south with probability $(1-p)/2$, respectively, with $p$ {being chosen in} $[0,1]$. Intuitively, random dispersal strategies, \textit{i.e.} $p=1/2$, 
might not be evolutionarily stable as the distribution of resources {is} generally not the same in {the} horizontal and vertical directions {so that it could} be more advantageous for {the} population to have {a} higher probability moving in one direction than the other. This is indeed the case, and {one might attempt to conjecture} that some mixed strategy $p^*$ {in} $(0,1)$, {allowing the} population to move in horizontal and vertical directions with different probabilities, {would} emerge as {an} evolutionarily stable dispersal strategy in {this} particular setting. A bit surprisingly, our results suggest that the only  evolutionarily stable dispersal strategies are $p^*=0$ and/or $p^*=1$, \textit{i.e.} it is  more advantageous for {the} population to move {in only} one direction.

More specifically, we consider a system of two competing populations in two-dimensio\-nal heterogeneous environments. The populations are assumed to move horizontally and vertically with different probabilities, but are otherwise identical. We introduce a function $F$ of the dispersal probability, which measures the difference between the spatial variations of the population distributions at equilibrium in horizontal and vertical directions: when it is positive, the species has more  variations in the horizontal direction; when it 
is negative, it has more variations in the vertical direction. 
We show that $F$ is monotone deceasing and that the  evolutionarily stable dispersal strategies are to maximize the function $F$ when it is positive and to minimize it when it is negative, {\it i.e.} the  evolutionarily stable strategies are to move in one direction. 
As the population distribution  is often positively correlated with the resource distribution, 
thus function $F$ indirectly measures the difference between the resource variations  in
horizontal and vertical directions. Therefore, our results seem to predict that it is more favorable for the species to choose the direction with smaller variation in the resource distribution.

To explain these findings intuitively, consider {a} peculiar scenario {in which} the underlying habitat is a rectangular region and the resources {are} distributed inhomogeneous{ly} in {the} horizontal direction but homogeneous{ly} in {the} vertical {one}. For such case, as there is only spatial variations in {the} horizontal direction, the results of Hastings \cite{Hastings1983} and Dockery \textit{et al.} \cite{DHMP} for the evolution of slow dispersal 
suggest that it might be better for the population not to move horizontally, 
which is in agreement with our findings in this paper. 
These considerations  also suggest that we should probably consider a broader {set} of dispersal strategies, \textit{e.g.} those strategies  which {allow for condition-dependent} movement \cite{CCL2010, Cosner2014, Gyllenberg2016, Kisdi, LamLou2014-1, LamLou2014-2, Potapov2014}.

\subsection{Organization of the paper}
Section \ref{sec:model} contains the formal derivations of the mathematical models at stage and the main results. 
We  present  numerical {simulations} in Section \ref{sec:numerics}
to complement the analytical results and to provide some intuition and insights. 
In Section \ref{sec:invasion fitness}, we discuss the stability of semi-trivial equilibria and investigate  properties of the invasion fitness. In Section \ref{sec:local stability}, we further study the  stability of the semi-trivial equilibria and 
identify all evolutionarily stable strategies. Section \ref{sec:full dynamics} is devoted to the classification of the global dynamics of the two-species competition model introduced in Section \ref{sec:model}. In Section \ref{sec:discussions}, we summarize our conclusions, 
and discuss possible extensions {of the present work}. {Finally}, some technical materials are given in the Appendix.

\section{The models and analytical results}\label{sec:model}

\subsection{Formal derivation from random walks}\label{subsec:formalder}
In this section, we will closely follow the approach in \cite{Slover}. Let the habitat be the discrete lattices of steps $\Delta x$ and $\Delta y$ in the full two-dimensional space $\R^2$. Assume that each individual moves horizontally to the left and to the right with probability $\tfrac{\theta}2$ and vertically up and down with probability $\tfrac{1-\theta}2$, with $\theta$ in $(0,1)$. Let $N(t,x,y)$ denote the number of individuals of the population at time $t$ and location $(x,y)$ in $\R^2$. Then
$$
\aligned
N(t+\Delta t,x,y)=&\frac{\theta}{2}\left[N(t,x+\Delta x,y)+N(t,x-\Delta x,y)\right] \\
&+\frac{1-\theta}{2}\left[N(t,x+\Delta x,y)+N(t,x-\Delta x,y)\right].
\endaligned
$$
Using Taylor series expansions, we have
\[
\aligned
\frac{N(t+\Delta t, x, y)-N(t, x, y)}{\Delta t}=&\frac{(\Delta x)^2}{2\Delta t} \theta N_{xx}+\frac{(\Delta y)^2}{2\Delta t} (1-\theta) N_{yy}\\
&+\frac{(\Delta x)^2}{\Delta t}\cdot O(\Delta x)+\frac{(\Delta y)^2}{\Delta t}\cdot O(\Delta y),
\endaligned
\]
where $O(\Delta x)$ and $O(\Delta y)$ denote terms which are bounded with respect to the $\Delta x$ and $\Delta y$, respectively. Assuming that both $\frac{(\Delta x)^2}{2\Delta t}\to \overline{D}$ and $\frac{(\Delta y)^2}{2\Delta t}\to\overline{D}$ for some positive constant $\overline{D}$ as the lengths of the time step $\Delta t$ and of the space steps $\Delta x$ and $\Delta y$ tend to $0$, we obtain, passing to the limit in the relation above,
\begin{equation}\label{eq:diffusiononly}
N_t=\overline{D}\left[\theta N_{xx}+(1-\theta)N_{yy}\right],
\end{equation}
which is the type of anisotropic {diffusion} operator to be considered in {the present} article.

The parameter $\theta$ can be regarded as a dispersal strategy. Namely, when $\theta=0$, the whole population will either move north or south with probability $1/2$. Similarly, when $\theta=1$, the population will only move east {or} west. Most of previous studies assume that individuals are randomly diffusing, \textit{i.e.} they move in two orthogonal directions with the same probability ($\theta=1/2$). Given arbitrarily distributed resources, is there some particular strategy $\theta$ {in} $[0,1]$ which can convey {a} competitive advantage? The main goal of {the present} article is to address this question.

\subsection{The single-species model}
{Incorporating} the population dynamics {into equation} \eqref{eq:diffusiononly}, we arrive at the following reaction-diffusion equation:
\begin{equation}\label{eq:singlespecies}
\left\{
\aligned
&N_t=D(\theta) N_{xx}+D(1-\theta)N_{yy}+(a-N)N\text{ {in} }\Omega,\ t>0,\\
&D(\theta)N_{x}\nu_x+D(1-\theta)N_{y}\nu_y=0\text{ {on} }\partial\Omega,\ t>0,\\
&N({0,\cdot,\cdot})=N_0\gneqq0\text{ {in} }\Omega.
\endaligned
\right.
\end{equation}

Here, the domain $\Omega$ is a smooth open bounded subset of $\R^2$, and we denote its boundary by $\partial \Omega$. We assume without any further notice that $\overline\Omega$ is strictly convex and $\partial\Omega$ {is of class} $C^1$. The unit outward normal vector to $\partial \Omega$ is denoted by $\nu=(\nu_x, \nu_y)$. Thus, the map $\nu$ is one-to-one and continuous from $\partial \Omega$ to $\mathbb{S}^1$.

For any $\theta$ {in} $[0,1]$, the function $D$ is defined by
\[
D(\theta):=\underline{D}+(\overline{D}-\underline{D})\theta,
\qquad \theta\in [0, 1].
\]
Note that if $\underline{D}=0$, $D(\theta)=\overline{D}\theta$ is reduced to the form
of diffusion in {equation} \eqref{eq:diffusiononly}, which is degenerate when $\theta=0$. To avoid such degeneracy, we assume in the remainder of the paper that $\overline{D}$ and $\underline{D}$ are positive constants satisfying
\[
0<\underline{D}<\overline{D}.
\]

The equation has been completed with zero flux 
boundary conditions so that no individuals may escape the domain. The model is not mass conservative since individuals may reproduce according to monostable non-linearities. The free growth rate is given by the heterogeneous function $a$, that is assumed to satisfy {the following assumption}:

\bigskip

\noindent{(A1)}
Function $(x,y) \mapsto a(x,y)$ is positive, H\"older continuous and non-constant in $\overline\Omega$.
 \bigskip


By standard regularity theory for parabolic equations and a comparison argument, see, \textit{e.g.} \cite{CC2003}, it can be shown that $N$ is positive {in} $\overline\Omega$ for all times and that $N(t,\cdot,\cdot)$ {tends to} $N_\theta$ uniformly {in} $\overline\Omega$ as $t$ {tends to infinity}, where $N_\theta$, {the equilibrium distribution of the population}, is the unique positive steady state of system \eqref{eq:singlespecies}, i.e. $N_\theta$ satisfies
\begin{equation}\label{eq:singlespeciessteady}
\left\{
\aligned
&D(\theta) (N_{\theta})_{xx}+D(1-\theta)(N_{\theta})_{yy}+(a-N_\theta)N_\theta=0\text{ {in} }\Omega,\\
&(D(\theta)(N_{\theta})_{x}, D(1-\theta) (N_{\theta})_{y})\cdot \nu=0\text{ {on} }\partial\Omega.
\endaligned
\right.
\end{equation}
Note that when {the function} $a$ is non-constant, {so is} $N_\theta$. 

Clearly, the parameter $\theta$ has {a} strong {influence} on $N_\theta$.
As $\theta$ increases, 
the single species has more tendency to move horizontally than vertically, which may reduce the spatial variations of the population distributions in the horizontal direction 
and increase 
the variations  in the vertical direction.
In this connection, we
have the following result: 
\begin{theorem}\label{thm:F-1} Define, for {$\theta$ in $[0,1]$}, the function
\begin{equation}\label{eq:functionF}
F(\theta):=\int_\Omega \left[\left((N_\theta)_{x}\right)^2-\left((N_\theta)_{y}\right)^2\right]\, dx\,dy.
\end{equation}
 Then $F'(\theta)<0$
for $\theta\in (0, 1)$. In particular, $F$ is either strictly positive, strictly negative or sign-changing exactly once in $(0, 1)$.
\end{theorem}
The function $F$, which plays a critical role in later analysis,  can be regarded as a measurement of  the difference between the variations of 
the population distributions in horizontal and vertical
directions: when $F$ is  positive, we envision that 
the species at equilibrium has more spatial variations horizontally;
when $F$ is negative, it has more variations in the vertical direction.
Theorem \ref{thm:F-1} 
implies that 
as the species increases the horizontal  diffusion   and reduces the vertical diffusion, 
then it tends to have more  variations in
the vertical direction than the horizontal direction.

As the population distribution  is often positively correlated with the resource distribution, function $F$ can also be viewed as an indirect measurement of the difference between the resource variations in the horizontal and vertical directions. Numerical results on the shape of the function $F$ are presented in the next section (see Figure \ref{fig:F}).  

\subsection{The two-species competition model}
Given arbitrarily distributed resources across the habitat, we may regard {the} parameter $\theta$ as a dispersal strategy and ask whether there is some {value for} $\theta$ which is evolutionarily stable. To address this question, we now move to the situation where two populations are competing for the same resources but adopt different dispersal strategies. We {thus} consider the following reaction-diffusion system for two competing species:
\begin{equation}\label{eq:main}
\left\{
\aligned
&U_t=D(p)U_{xx}+D(1-p)U_{yy}+(a-U-V)U\text{ {in} }\Omega,\ t>0,\\
&V_t=D(q)V_{xx}+D(1-q)V_{yy}+(a-U-V)V\text{ {in} }\Omega,\ t>0,\\
&(D(p)U_x,D(1-p)U_y)\cdot\nu=0\text{ {on} }\partial\Omega,\ t>0,\\
&(D(q)V_x,D(1-q)V_y)\cdot\nu=0\text{ {on} }\partial\Omega,\ t>0,\\
&U(0,\cdot,\cdot)=U_0\gneqq 0,\ V(0,\cdot,\cdot)=V_0\gneqq 0\text{ {in} }\Omega,
\endaligned
\right.
\end{equation}
{in which the functions} $U$ and $V$ represent the {respective} population densities of two competing species. By standard regularity theory and {the} maximum principle for parabolic equations, it can be shown that $U$ and $V$ are positive {in} $\overline\Omega$ for all times. The competition for resources is neutral and independent of the dispersal strategy of the individuals, so that the death rate is given by $U+V$ for both populations.

We are given two orthogonal space directions $e_1=(1,0)$ and $e_2=(0,1)$, so that $(x,y)$ are the Cartesian coordinates in this basis. We may emphasise that after that choice, the problem is not rotationally invariant. As such, the two populations disperse with their own dispersal strategies, assimilated to the {respective probabilities} $p$ {and} $q$ to move in the direction $e_1$, with $p$ {and $q$ chosen in $[0,1]$}. As formally explained in {Subsection} \ref{subsec:formalder}, this way of dispersing results in a diffusion coefficient given by $D(p)$ (resp. $D(q)$) in the direction $e_1$ and $D(1-p)$ (resp. $D(1-q)$) in the direction $e_2$ {for the first (resp. second) density}.




We will adopt the viewpoint in the theory of adaptive dynamics. An important concept in adaptive dynamics is that of \emph{evolutionarily stable strategies} (ESS). A strategy is said to be evolutionarily stable if a population using it cannot be invaded by any small population using a different strategy. In {system} \eqref{eq:main}, $p$ and $q$ are strategies for two populations. In terms of adaptive dynamics, we say that  $p$ in $[0, 1]$ is an ESS if the semi-trivial steady state $(N_p,0)$ is locally asymptotically stable for $q\neq p$, with $q$ in $[0,1]$ and $q$ close to $p$. 

The following result characterizes the local stability of $(N_p 0)$ for $p$ and $q$ in $[0,1]$. 

\begin{theorem}\label{thm:Fpositive-001}
There exists some continuous function $q=q^*(p)$, defined in $[0,1]$, satisfying $0\le q^*(p)\le 1$ such that the following statements hold. 

\begin{enumerate}[\rm (i)]
\item If $F$ is positive in $[0,1]$, then $q^*(p)>p$ and $q^*(p)\equiv 1$ for $p$ close to $1$ such that $(N_p, 0)$ is stable for $p<q<q^*(p)$, unstable for $q>q^*(p)$ and $q<p$. In particular, if $q^*(p)\equiv 1$, then $(N_p,0)$ is stable for $q>p$, unstable for $q<p$.\label{thm:Fpositive-001 - statement positive}

\medskip

\item If $F$ has a unique root $\theta^*$ in $(0,1)$, then $q^*(p)>p$ for $p\in [0,\theta^*)$ and $q^*(p)<p$ for $p\in (\theta^*,1]$, such that $(N_p, 0)$ is stable for  $\min\{q^*(p), p\}<q<\max\{q^*(p), p\}$, unstable for $q>\max\{q^*(p), p\}$ and $q<\min\{q^*(p),p\}$.\label{thm:Fpositive-001 - statement root}

\medskip
\item If $F$ is negative in $[0, 1]$, then $q^*(p)<p$ and $q^*(p)\equiv 0$ for $p$ close to $0$ such that $(N_p, 0)$ is stable for $q^*(p)<q<p$, and unstable  for $q<q^*(p)$ and  $q>p$. In particular, if $q^*(p)\equiv 0$, then $(N_p, 0)$ is stable for $q<p$ and unstable for $q>p$.\label{thm:Fpositive-001 - statement negative}
\end{enumerate} 
\end{theorem}

This result follows from Theorems \ref{thm:Fpositive-1}, \ref{thm:Fpositive-2}, and \ref{thm:Fpositive-3}. In the next section, numerical results shed some insight into the stability of $(N_p,0)$ and illustrate the conclusions of Theorem \ref{thm:Fpositive-001} (see Figure \ref{fig:nodal set single}). Some biological intuition can also be gained from this Theorem as it provides a criterion for finding the ESS of system~\eqref{eq:main}.

\begin{corollary}\label{thm:F}
The following conclusions hold.

\begin{enumerate}[\noindent\rm (i)]
\item If $F$ is positive in $[0,1]$, then $p=0$ is the only ESS.\label{thm:F - statement positive}
\medskip

\item If $F$ has exactly one root $\theta^*$ in $(0, 1)$ so that $F$ is positive in $[0,\theta^*)$ and negative in $(\theta^*, 1]$, {then} both $p=0$ and $p=1$ are ESS, and $\theta^*$ is not evolutionarily stable.\label{thm:F - statement root}
\medskip

\item If {$F$ is negative} in $[0,1]$, and $p=1$ is the only ESS.\label{thm:F - statement negative}
\end{enumerate}
\end{corollary}

Our remaining goals include understanding the global dynamics of system \eqref{eq:main}. {This} system {possesses} two semi-trivial steady states, given by $(N_p,0)$ and $(0,N_q)$, respectively. Theorem \ref{thm:Fpositive-001} addresses the local stability of $(N_p,0)$ for arbitrary values of $p$ and $q$, and the stability of $(0, N_q)$ can be similarly determined. Furthermore, we shall show that there are only three alternatives for the global dynamics of system \eqref{eq:main}:
\begin{enumerate}[(i)]
\item {the state} $(N_p,0)$ is globally stable;
\item {the state} $(0, N_q)$ is globally stable;
\item {the states} $(N_p,0)$ and $(0, N_q)$ are both unstable, and there exists a unique positive steady state which is globally stable.
\end{enumerate}

We refer to the statements of Theorems \ref{thm:Fpositive-1-aa}, \ref{thm:Fpositive-2-aa} and \ref{thm:Fchangessign-1} for further details on the characterizations of the global dynamics of system \eqref{eq:main}. These analytical results on the dynamics of \eqref{eq:main}, complemented by numerical simulations in the next section for {a free growth rate function of the form} $a(x,y)=\lambda A(x)+(1-\lambda)A(y)$, with $\lambda$ {in} $[0,1]$, will help provide a more clear picture on the dynamics of system \eqref{eq:main}.

\section{The numerical results}\label{sec:numerics}

All the simulations presented here 
were achieved using the free and open-source software \textsc{FreeFEM} \cite{Hecht:2012}. The numerical approximation of the large-time solution to system \eqref{eq:singlespecies} was based on a variational form of the problem and achieved using a spatial discretisation based on the finite element method, with $P_1$ Lagrange elements, combined with an implicit-explicit (IMEX) Euler scheme (see \cite{Asher:1995} for instance) for the time-integration of the resulting ordinary differential equations. The linear terms in the reaction-diffusion equation are then treated implicitly in time, while the non-linear reaction term is dealt with explicitly, in order to enforce the stability of the scheme.

The mesh used to discretise the domain $\Omega$ realised as a disk of radius $2$ was comprised of $3916$ triangles and the length of the time step used was $0.01$. The chosen initial state $N_0$ is the constant one, with value $0.5$. Once the stationarity of the approximate solution was obtained in relative $L_2$ norm within a prescribed tolerance of $1\mathrm{e}{-15}$, the approximate steady state was used to compute an approximate value of $F(p)$ and also for a finite element discretisation of the linear eigenvalue problem \eqref{eq:varphi-1}. The computation of an approximation of the smallest eigenvalue of \eqref{eq:varphi-1}, denoted by $\Lambda(p,q)$, was done with the ARPACK package. Representations of an approximation to the nodal set of $\Lambda(p,q)$ for different values of $\lambda$ were then obtained by repeating the computation for numerous values of the parameters $p$ and $q$ taken in a discrete grid of the interval $[0,1]$.

Note that the state $(N_p,0)$ is linearly stable when $\Lambda(p,q)$ is positive, and unstable when $\Lambda(p,q)$ is negative. Furthermore, $\Lambda(p,q)$ vanishes whenever $p=q$, \textit{i.e.} the nodal sets of $\Lambda(p,q)$ always consist of the diagonal line $p=q$ in the $p-q$ plane. 

\subsection{The function $F$}
Numerical approximations of the graph of the function $F$ in the case of a free-growth function of the form $a(x,y)=\lambda A(x)+(1-\lambda)A(y)$ are provided in Figure \ref{fig:F}. For the simulations, we considered a disk of radius $2$ centered at the origin for the domain $\Omega$, anisotropic diffusion parameters $\underline{D}$ and $\overline{D}$ respectively equal to $0.1$ and $10$, and the function $A(x)=2-\sin(\pi x)$.

It is easily seen that, for all $\lambda$ in $[0,1]$ and all $\theta$ in $[0,1]$, the value $F(\theta)$ for $\lambda$ is equal to the value of $-F(1-\theta)$ for $1-\lambda$. Due to this symmetry in the function $F$ with respect to $\lambda$, we only plot the graph of $F$ for different values of $\lambda$ between $0$ and $0.5$, illustrating how the function goes from strictly negative, to sign-changing once, and to strictly positive as $\lambda$ varies.
As shown in Figure~\ref{fig:F}, the function $F$ is strictly decreasing in $\theta$, as predicted by Theorem~\ref{thm:F-1}.

\begin{figure}
\centering
\begin{subfigure}{0.45\textwidth}
\resizebox{\textwidth}{!}{\begin{tikzpicture}
\begin{axis}[xmin=0,xmax=1,xlabel=$\theta$,ylabel=$F(\theta)$]
\addplot[] table {case1_l0F.txt} ;
\end{axis}
\end{tikzpicture}}
\caption{$\lambda=0$\label{fig:F - lambda=0}}
\end{subfigure}
\begin{subfigure}{0.45\textwidth}
\resizebox{\textwidth}{!}{\begin{tikzpicture}
\begin{axis}[xmin=0,xmax=1,xlabel=$\theta$,ylabel=$F(\theta)$]
\addplot[] table {case1_l2F.txt} ;
\end{axis}
\end{tikzpicture}}
\caption{$\lambda=0.1$\label{fig:F - lambda=0.1}}
\end{subfigure}

\begin{subfigure}{0.45\textwidth}
\resizebox{\textwidth}{!}{\begin{tikzpicture}
\begin{axis}[xmin=0,xmax=1,xlabel=$\theta$,ylabel=$F(\theta)$]
\addplot[] table {case1_l4F.txt} ;
\end{axis}
\end{tikzpicture}}
\caption{$\lambda=0.2$}
\end{subfigure}
\begin{subfigure}{0.45\textwidth}
\resizebox{\textwidth}{!}{\begin{tikzpicture}
\begin{axis}[xmin=0,xmax=1,xlabel=$\theta$,ylabel=$F(\theta)$]
\addplot[] table {case1_l6F.txt} ;
\end{axis}
\end{tikzpicture}}
\caption{$\lambda=0.3$}
\end{subfigure}

\begin{subfigure}{0.45\textwidth}
\resizebox{\textwidth}{!}{\begin{tikzpicture}
\begin{axis}[xmin=0,xmax=1,xlabel=$\theta$,ylabel=$F(\theta)$]
\addplot[] table {case1_l8F.txt} ;
\end{axis}
\end{tikzpicture}}
\caption{$\lambda=0.4$}
\end{subfigure}
\begin{subfigure}{0.45\textwidth}
\resizebox{\textwidth}{!}{\begin{tikzpicture}
\begin{axis}[xmin=0,xmax=1,xlabel=$\theta$,ylabel=$F(\theta)$]
\addplot[] table {case1_l10F.txt} ;
\end{axis}
\end{tikzpicture}}
\caption{$\lambda=0.5$\label{fig:F - lambda=0.5}}
\end{subfigure}
\caption{Numerical approximations of the graphs of the function $F$ for $\lambda$ taking the values $0$, $0.1$, $0.2$, $0.3$, $0.4$, and $0.5$, in the case where $\Omega$ is a disk of radius $2$ centered at the origin, $\underline{D}=0.1$, $\overline{D}=10$, and $A(x)=2-\sin(\pi x)$.}\label{fig:F}
\end{figure}
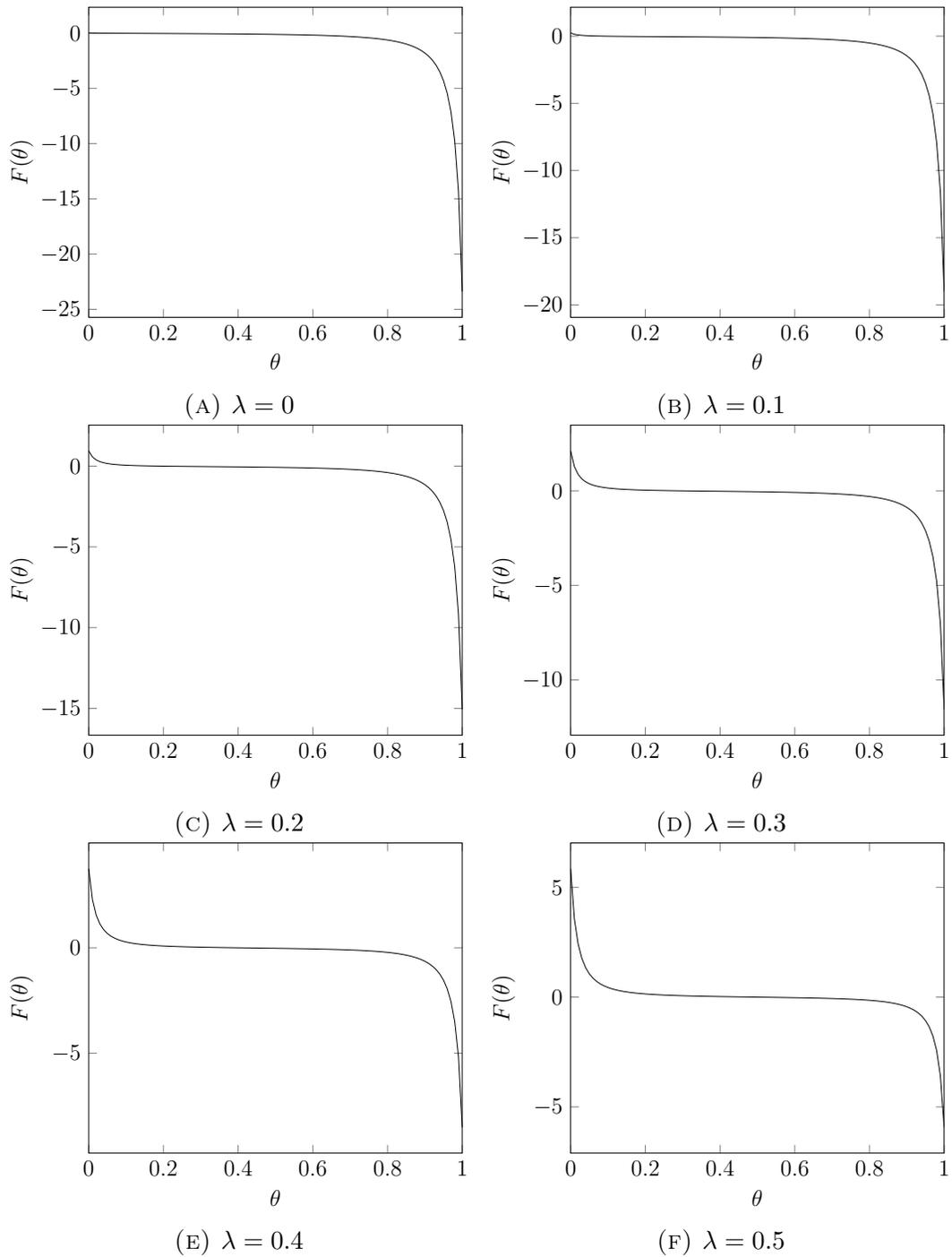

\medskip

Figures \ref{steady states, lambda=0.4} and \ref{steady states, lambda=0.6} present the numerical approximations of the function $a$ and of the steady state $N_\theta$ for various values of $\theta$ in $[0,1]$ and $\lambda$ respectively equal to $0.4$ and $0.6$, in the case where $\Omega$ is a disk of radius $2$ centered at the origin, $\underline{D}=0.1$, $\overline{D}=10$, and $A(x)=2-\sin(\pi x)$. For such values of $\lambda$, the function does not appear more biased in the horizontal direction than in the vertical one. Yet, one can see clearly that the steady state shows very little variation in vertical direction for $\theta=0$, but changes as the value of the parameter $\theta$ grows to end with very little variation in the horizontal direction for $\theta=1$. This illustrates how the function $F$ goes from being negative to positive as $\theta$ varies.

\begin{figure}
\centering
\begin{subfigure}{0.42\textwidth}
\resizebox{\textwidth}{!}{\input{case1_function_a_l8_seismic.tex}}
\caption{The function $a$.}
\end{subfigure}
	
\begin{subfigure}{0.42\textwidth}
\resizebox{\textwidth}{!}{\input{case1_steady_state_l8_p0_seismic.tex}}
\caption{The steady state for $\theta=0$.}
\end{subfigure}
\begin{subfigure}{0.42\textwidth}
\resizebox{\textwidth}{!}{\input{case1_steady_state_l8_p2_seismic.tex}}
\caption{The steady state for $\theta=0.2$.}
\end{subfigure}
	
\begin{subfigure}{0.42\textwidth}
\resizebox{\textwidth}{!}{\input{case1_steady_state_l8_p4_seismic.tex}}
\caption{The steady state for $\theta=0.4$.}
\end{subfigure}
\begin{subfigure}{0.42\textwidth}
\resizebox{\textwidth}{!}{\input{case1_steady_state_l8_p6_seismic.tex}}
\caption{The steady state for $\theta=0.6$.}
\end{subfigure}
	
\begin{subfigure}{0.42\textwidth}
\resizebox{\textwidth}{!}{\input{case1_steady_state_l8_p8_seismic.tex}}
\caption{The steady state for $\theta=0.8$.}
\end{subfigure}
\begin{subfigure}{0.42\textwidth}
\resizebox{\textwidth}{!}{\input{case1_steady_state_l8_p10_seismic.tex}}
\caption{The steady state for $\theta=1$.}
\end{subfigure}
\caption{Isolines for numerical approximations of the function $a$ and of the steady state $N_\theta$ for $\theta$ taking the values $0$, $0.2$, $0.4$, $0.6$, $0.8$, and $1$, in the case where $\Omega$ is a disk of radius $2$ centered at the origin, $\underline{D}=0.1$, $\overline{D}=10$, $\lambda=0.4$, and $A(x)=2-\sin(\pi x)$.}\label{steady states, lambda=0.4}
\end{figure}

\begin{figure}
\centering
\begin{subfigure}{0.42\textwidth}
\resizebox{\textwidth}{!}{\input{case1_function_a_l12_seismic.tex}}
\caption{The function $a$.}
\end{subfigure}
	
\begin{subfigure}{0.42\textwidth}
\resizebox{\textwidth}{!}{\input{case1_steady_state_l12_p0_seismic.tex}}
\caption{The steady state for $\theta=0$.}
\end{subfigure}
\begin{subfigure}{0.42\textwidth}
\resizebox{\textwidth}{!}{\input{case1_steady_state_l12_p2_seismic.tex}}
\caption{The steady state for $\theta=0.2$.}
\end{subfigure}
	
\begin{subfigure}{0.42\textwidth}
\resizebox{\textwidth}{!}{\input{case1_steady_state_l12_p4_seismic.tex}}
\caption{The steady state for $\theta=0.4$.}
\end{subfigure}
\begin{subfigure}{0.42\textwidth}
\resizebox{\textwidth}{!}{\input{case1_steady_state_l12_p6_seismic.tex}}
\caption{The steady state for $\theta=0.6$.}
\end{subfigure}
	
\begin{subfigure}{0.42\textwidth}
\resizebox{\textwidth}{!}{\input{case1_steady_state_l12_p8_seismic.tex}}
\caption{The steady state for $\theta=0.8$.}
\end{subfigure}
\begin{subfigure}{0.42\textwidth}
\resizebox{\textwidth}{!}{\input{case1_steady_state_l12_p10_seismic.tex}}
\caption{The steady state for $\theta=1$.}
\end{subfigure}
\caption{Isolines for numerical approximations of the function $a$ and of the steady state $N_\theta$ for $\theta$ taking the values $0$, $0.2$, $0.4$, $0.6$, $0.8$, and $1$, in the case where $\Omega$ is a disk of radius $2$ centered at the origin, $\underline{D}=0.1$, $\overline{D}=10$, $\lambda=0.6$, and $A(x)=2-\sin(\pi x)$.}\label{steady states, lambda=0.6}
\end{figure}

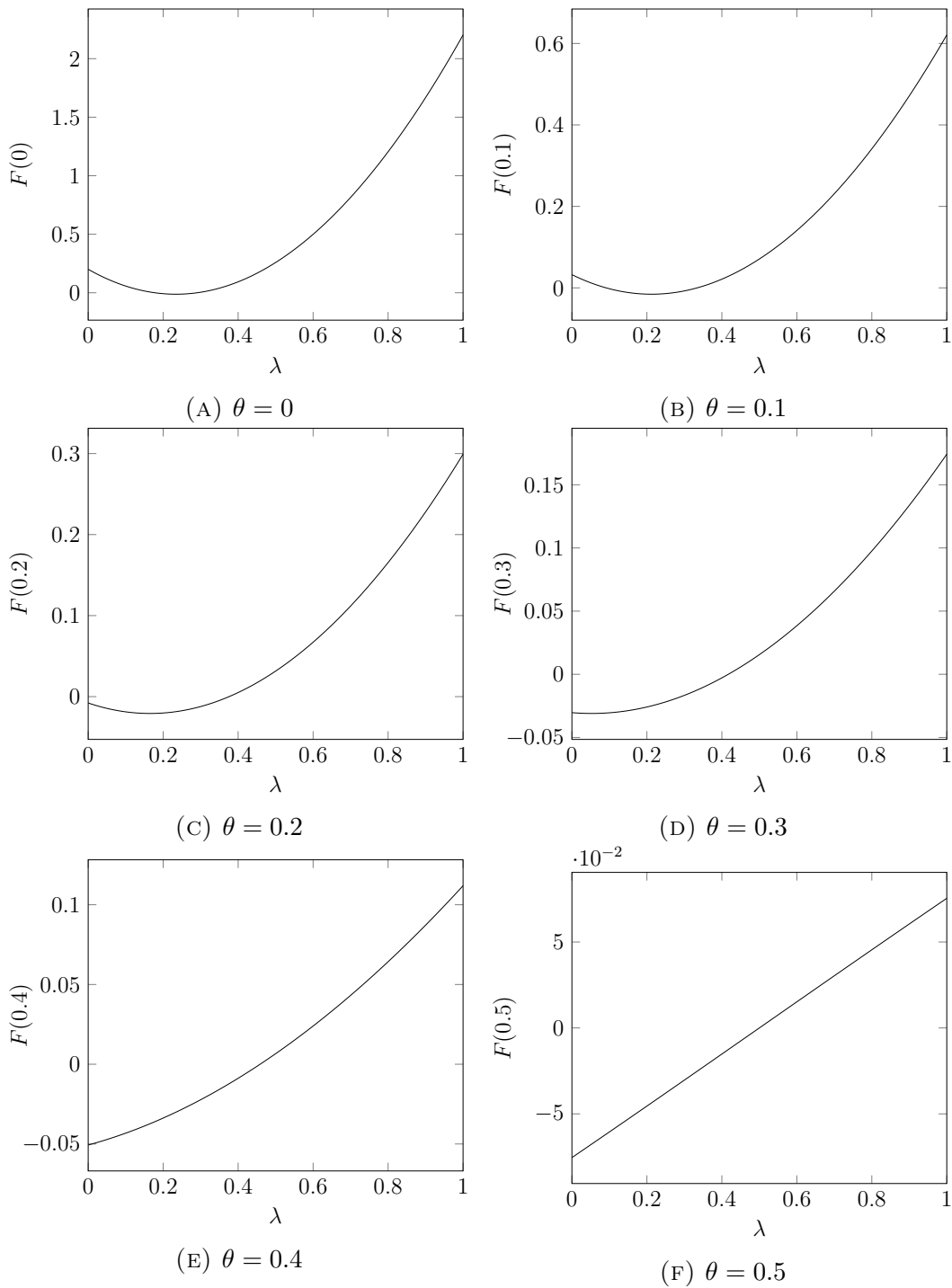
\begin{figure}
\centering
\begin{subfigure}{0.45\textwidth}
\resizebox{\textwidth}{!}{\begin{tikzpicture}
\begin{axis}[xmin=0,xmax=1,xlabel=$\lambda$,ylabel=$F(0)$,tick label style={/pgf/number format/fixed}]
\addplot[] table {case2_Fp00.txt} ;
\end{axis}
\end{tikzpicture}}
\caption{$\theta=0$}
\end{subfigure}
\begin{subfigure}{0.45\textwidth}
\resizebox{\textwidth}{!}{\begin{tikzpicture}
\begin{axis}[xmin=0,xmax=1,xlabel=$\lambda$,ylabel=$F(0.1)$,tick label style={/pgf/number format/fixed}]
\addplot[] table {case2_Fp01.txt} ;
\end{axis}
\end{tikzpicture}}
\caption{$\theta=0.1$}
\end{subfigure}

\begin{subfigure}{0.45\textwidth}
\resizebox{\textwidth}{!}{\begin{tikzpicture}
\begin{axis}[xmin=0,xmax=1,xlabel=$\lambda$,ylabel=$F(0.2)$,tick label style={/pgf/number format/fixed}]
\addplot[] table {case2_Fp02.txt} ;
\end{axis}
\end{tikzpicture}}
\caption{$\theta=0.2$}
\end{subfigure}
\begin{subfigure}{0.45\textwidth}
\resizebox{\textwidth}{!}{\begin{tikzpicture}
\begin{axis}[xmin=0,xmax=1,xlabel=$\lambda$,ylabel=$F(0.3)$,tick label style={/pgf/number format/fixed}]
\addplot[] table {case2_Fp03.txt} ;
\end{axis}
\end{tikzpicture}}
\caption{$\theta=0.3$}
\end{subfigure}

\begin{subfigure}{0.45\textwidth}
\resizebox{\textwidth}{!}{\begin{tikzpicture}
\begin{axis}[xmin=0,xmax=1,xlabel=$\lambda$,ylabel=$F(0.4)$,tick label style={/pgf/number format/fixed}]
\addplot[] table {case2_Fp04.txt} ;
\end{axis}
\end{tikzpicture}}
\caption{$\theta=0.4$}
\end{subfigure}
\begin{subfigure}{0.45\textwidth}
\resizebox{\textwidth}{!}{\begin{tikzpicture}
\begin{axis}[xmin=0,xmax=1,xlabel=$\lambda$,ylabel=$F(0.5)$,tick label style={/pgf/number format/fixed}]
\addplot[] table {case2_Fp05.txt} ;
\end{axis}
\end{tikzpicture}}
\caption{$\theta=0.5$}
\end{subfigure}
\caption{Numerical approximations of the graphs of $F(\theta)$ as a function of $\lambda$ for $\theta$ taking the values $0$, $0.1$, $0.2$, $0.3$, $0.4$, and $0.5$, in the case where $\Omega$ is a disk of radius $2$ centered at the origin, $\underline{D}=0.1$, $\overline{D}=10$, and $A(x)=4-\frac{1}{4}\,x^2$.}\label{fig:Fp}
\end{figure}

In Figure \ref{fig:Fp}, we fix a value of $\theta$ and observe that $F(\theta)$, seen as a function of $\lambda$, is not necessarily monotone. In this case, we considered a disk of radius $2$ centered at the origin for the domain $\Omega$, anisotropic diffusion parameters $\underline{D}$ and $\overline{D}$ respectively equal to $0.1$ and $10$, and the function $A(x)=4-\frac{1}{4}\,x^2$.

\subsection{Local stability of $(N_p,0)$}
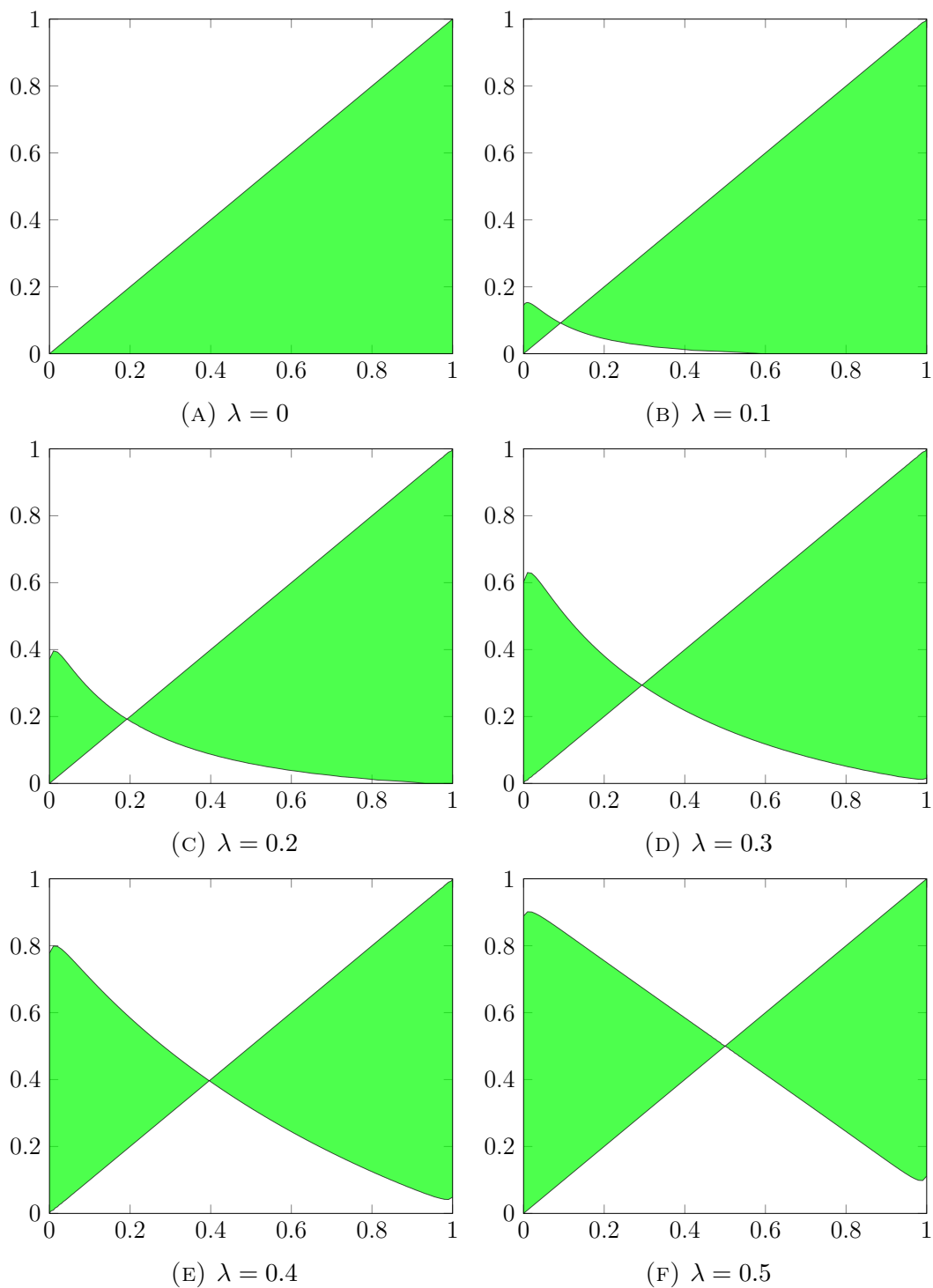
\begin{figure}
\centering
\begin{subfigure}{0.45\textwidth}
\resizebox{\textwidth}{!}{\begin{tikzpicture}
\begin{axis}[xmin=0,xmax=1,ymin=0,ymax=1]
\addplot[fill=green,opacity=0.7] table {case1_l0PQ1.txt} --cycle;
\end{axis}
\end{tikzpicture}}
\caption{$\lambda=0$\label{fig:nodal set single - lambda=0}}
\end{subfigure}
\begin{subfigure}{0.45\textwidth}
\resizebox{\textwidth}{!}{\begin{tikzpicture}
\begin{axis}[xmin=0,xmax=1,ymin=0,ymax=1]
\addplot[fill=green,opacity=0.7] table {case1_l2PQ1.txt} --cycle;
\addplot[fill=green,opacity=0.7] table {case1_l2PQ2.txt} --cycle;
\end{axis}
\end{tikzpicture}}
\caption{$\lambda=0.1$\label{fig:nodal set single - lambda=0.1}}
\end{subfigure}

\begin{subfigure}{0.45\textwidth}
\resizebox{\textwidth}{!}{\begin{tikzpicture}
\begin{axis}[xmin=0,xmax=1,ymin=0,ymax=1]
\addplot[fill=green,opacity=0.7] table {case1_l4PQ1.txt} --cycle;
\addplot[fill=green,opacity=0.7] table {case1_l4PQ2.txt} --cycle;
\end{axis}
\end{tikzpicture}}
\caption{$\lambda=0.2$}
\end{subfigure}
\begin{subfigure}{0.45\textwidth}
\resizebox{\textwidth}{!}{\begin{tikzpicture}
\begin{axis}[xmin=0,xmax=1,ymin=0,ymax=1]
\addplot[fill=green,opacity=0.7] table {case1_l6PQ1.txt} --cycle;
\addplot[fill=green,opacity=0.7] table {case1_l6PQ2.txt} --cycle;
\end{axis}
\end{tikzpicture}}
\caption{$\lambda=0.3$}
\end{subfigure}

\begin{subfigure}{0.45\textwidth}
\resizebox{\textwidth}{!}{\begin{tikzpicture}
\begin{axis}[xmin=0,xmax=1,ymin=0,ymax=1]
\addplot[fill=green,opacity=0.7] table {case1_l8PQ1.txt} --cycle;
\addplot[fill=green,opacity=0.7] table {case1_l8PQ2.txt} --cycle;
\end{axis}
\end{tikzpicture}}
\caption{$\lambda=0.4$}
\end{subfigure}
\begin{subfigure}{0.45\textwidth}
\resizebox{\textwidth}{!}{\begin{tikzpicture}
\begin{axis}[xmin=0,xmax=1,ymin=0,ymax=1]
\addplot[fill=green,opacity=0.7] table {case1_l10PQ1.txt} --cycle;
\addplot[fill=green,opacity=0.7] table {case1_l10PQ2.txt} --cycle;
\end{axis}
\end{tikzpicture}}
\caption{$\lambda=0.5$\label{fig:nodal set single - lambda=0.5}}
\end{subfigure}
\caption{Numerical approximations of the nodal sets of $\Lambda(p,q)$ for $\lambda$ taking the values $0$, $0.1$, $0.2$, $0.3$, $0.4$, and $0.5$, in the case where $\Omega$ is a disk of radius $2$ centered at the origin, $\underline{D}=0.1$, $\overline{D}=10$, and $A(x)=2-\sin(\pi x)$. The subset colored in green is the one in which $\Lambda(p,q)>0$, that is in which the steady state $(N_p,0)$ is linearly stable.}\label{fig:nodal set single}
\end{figure}

Figure \ref{fig:nodal set single} presents numerical approximations of the nodal sets of $\Lambda(p,q)$ for a free growth function of the form $a(x,y)=\lambda A(x)+(1-\lambda)A(y)$, the choice of the problem parameters being the same as for the graphs of $F$ in Figure~\ref{fig:F}. 

Here, $\Lambda(p,q)$ is positive if and only if $(N_p,0)$ is linearly stable and  $\Lambda(p,q)<0$ if and only if $(0, N_q)$  is linearly unstable. Again, for all $\lambda$ in $[0,1]$ and all $(p,q)$ in $[0,1]^2$, the value of $\Lambda(p,q)$ for $\lambda$ is equal to the value of $\Lambda(1-p,1-q)$ for $1-\lambda$. Due to this symmetry of $\Lambda(p,q)$ with respect to $\lambda$, we only plot the nodal sets for values of $\lambda$ between $0$ and $0.5$, providing a graphical illustration of how the nodal set of $\Lambda(p,q)$ changes as $\lambda$ varies, in connection with the results of Theorem \ref{thm:Fpositive-001}, or more broadly, the conclusions of Theorems \ref{thm:Fpositive-1}, \ref{thm:Fpositive-2}, and~\ref{thm:Fpositive-3}.

Subfigure \ref{fig:nodal set single}\subref{fig:nodal set single - lambda=0} corresponds to statement \eqref{thm:Fpositive-001 - statement negative} of Theorem \ref{thm:Fpositive-001} with $q^*(p)\equiv 0$, which shows that $(N_p,0)$ is stable for $p>q$ and unstable for $p<q$. In particular, $p=1$ is the only ESS. In this case, the nodal set of $\Lambda(p,q)$ consists precisely of the diagonal line $q=p$. These conclusions are consistent with Subfigure \ref{fig:F}\subref{fig:F - lambda=0}, in which $F$ is strictly negative. 

Subfigures \ref{fig:nodal set single}\subref{fig:nodal set single - lambda=0.1} to \ref{fig:nodal set single}\subref{fig:nodal set single - lambda=0.5} correspond to statement \eqref{thm:Fpositive-001 - statement root} of Theorem \ref{thm:Fpositive-001}, which shows that $(N_p,0)$ is stable for $\min\{q^*(p),p\}<q<\max\{q^*(p),p\}$, unstable for $\max\{q^*(p),p\}<q$ and $q<\min\{q^*(p),p\}$, and the nodal set of $\Lambda(p,q)$ consists of the curve $q=q^*(p)$ and the diagonal line $q=p$. Hence, both $p=0$ and $p=1$ are the only ESS. These conclusions are in accordance with Subfigures 
\ref{fig:F}\subref{fig:F - lambda=0.1} to 
\ref{fig:F}\subref{fig:F - lambda=0.5}, in which $F$ changes sign exactly once in $(0,1)$, from positive to negative.

\subsection{Global dynamics}
\begin{figure}
\centering
\begin{subfigure}{0.4\textwidth}
\resizebox{\textwidth}{!}{\begin{tikzpicture}
\begin{axis}[xmin=0,xmax=1,ymin=0,ymax=1]
\addplot[fill=green,opacity=0.7] table {case2_l0PQ1.txt} --cycle;
\addplot[fill=green,opacity=0.7] table {case2_l0PQ2.txt} --cycle;
\addplot[fill=red,opacity=0.7] table {case2_l0QP1.txt} --cycle;
\addplot[fill=red,opacity=0.7] table {case2_l0QP2.txt} --cycle;
\end{axis}
\end{tikzpicture}}
\caption{$\lambda=0$}
\end{subfigure}
\begin{subfigure}{0.4\textwidth}
\resizebox{\textwidth}{!}{\begin{tikzpicture}
\begin{axis}[xmin=0,xmax=1,ymin=0,ymax=1]
\addplot[fill=green,opacity=0.7] table {case2_l2PQ1.txt} --cycle;
\addplot[fill=green,opacity=0.7] table {case2_l2PQ2.txt} --cycle;
\addplot[fill=red,opacity=0.7] table {case2_l2QP1.txt} --cycle;
\addplot[fill=red,opacity=0.7] table {case2_l2QP2.txt} --cycle;
\end{axis}
\end{tikzpicture}}
\caption{$\lambda=0.1$}
\end{subfigure}

\begin{subfigure}{0.4\textwidth}
\resizebox{\textwidth}{!}{\begin{tikzpicture}
\begin{axis}[xmin=0,xmax=1,ymin=0,ymax=1]
\addplot[fill=green,opacity=0.7] table {case2_l4PQ1.txt} --cycle;
\addplot[fill=red,opacity=0.7] table {case2_l4QP1.txt} --cycle;
\end{axis}
\end{tikzpicture}}
\caption{$\lambda=0.2$}
\end{subfigure}
\begin{subfigure}{0.4\textwidth}
\resizebox{\textwidth}{!}{\begin{tikzpicture}
\begin{axis}[xmin=0,xmax=1,ymin=0,ymax=1]
\addplot[fill=green,opacity=0.7] table {case2_l6PQ1.txt} --cycle;
\addplot[fill=green,opacity=0.7] table {case2_l6PQ2.txt} --cycle;
\addplot[fill=red,opacity=0.7] table {case2_l6QP1.txt} --cycle;
\addplot[fill=red,opacity=0.7] table {case2_l6QP2.txt} --cycle;
\end{axis}
\end{tikzpicture}}
\caption{$\lambda=0.3$}
\end{subfigure}

\begin{subfigure}{0.4\textwidth}
\resizebox{\textwidth}{!}{\begin{tikzpicture}
\begin{axis}[xmin=0,xmax=1,ymin=0,ymax=1]
\addplot[fill=green,opacity=0.7] table {case2_l8PQ1.txt} --cycle;
\addplot[fill=green,opacity=0.7] table {case2_l8PQ2.txt} --cycle;
\addplot[fill=red,opacity=0.7] table {case2_l8QP1.txt} --cycle;
\addplot[fill=red,opacity=0.7] table {case2_l8QP2.txt} --cycle;
\end{axis}
\end{tikzpicture}}
\caption{$\lambda=0.4$}
\end{subfigure}
\begin{subfigure}{0.4\textwidth}
\resizebox{\textwidth}{!}{\begin{tikzpicture}
\begin{axis}[xmin=0,xmax=1,ymin=0,ymax=1]
\addplot[fill=green,opacity=0.7] table {case2_l10PQ1.txt} --cycle;
\addplot[fill=green,opacity=0.7] table {case2_l10PQ2.txt} --cycle;
\addplot[fill=red,opacity=0.7] table {case2_l10QP1.txt} --cycle;
\addplot[fill=red,opacity=0.7] table {case2_l10QP2.txt} --cycle;
\end{axis}
\end{tikzpicture}}
\caption{$\lambda=0.5$}
\end{subfigure}
\caption{Numerical approximations of the nodal sets of $\Lambda(p,q)$ and $\Lambda(q,p)$ for $\lambda$ taking the values $0$, $0.1$, $0.2$, $0.3$, $0.4$, and $0.5$, in the case where $\Omega$ is a disk of radius $2$ centered at the origin, $\underline{D}=0.1$, $\overline{D}=10$, and $A(x)=4-\frac{1}{4}\,x^2$. The subset colored in green is the one in which $\Lambda(p,q)>0>\Lambda(q,p)$, i.e. in which  $(N_p,0)$ is stable, while that colored in red is the one in which $\Lambda(q,p)>0>\Lambda(q,p)$, i.e. in which  $(0,N_q)$ is stable. The white region is where both $(N_p,0)$ and $(0,N_q)$ are unstable.}\label{fig: nodal set double}
\end{figure}
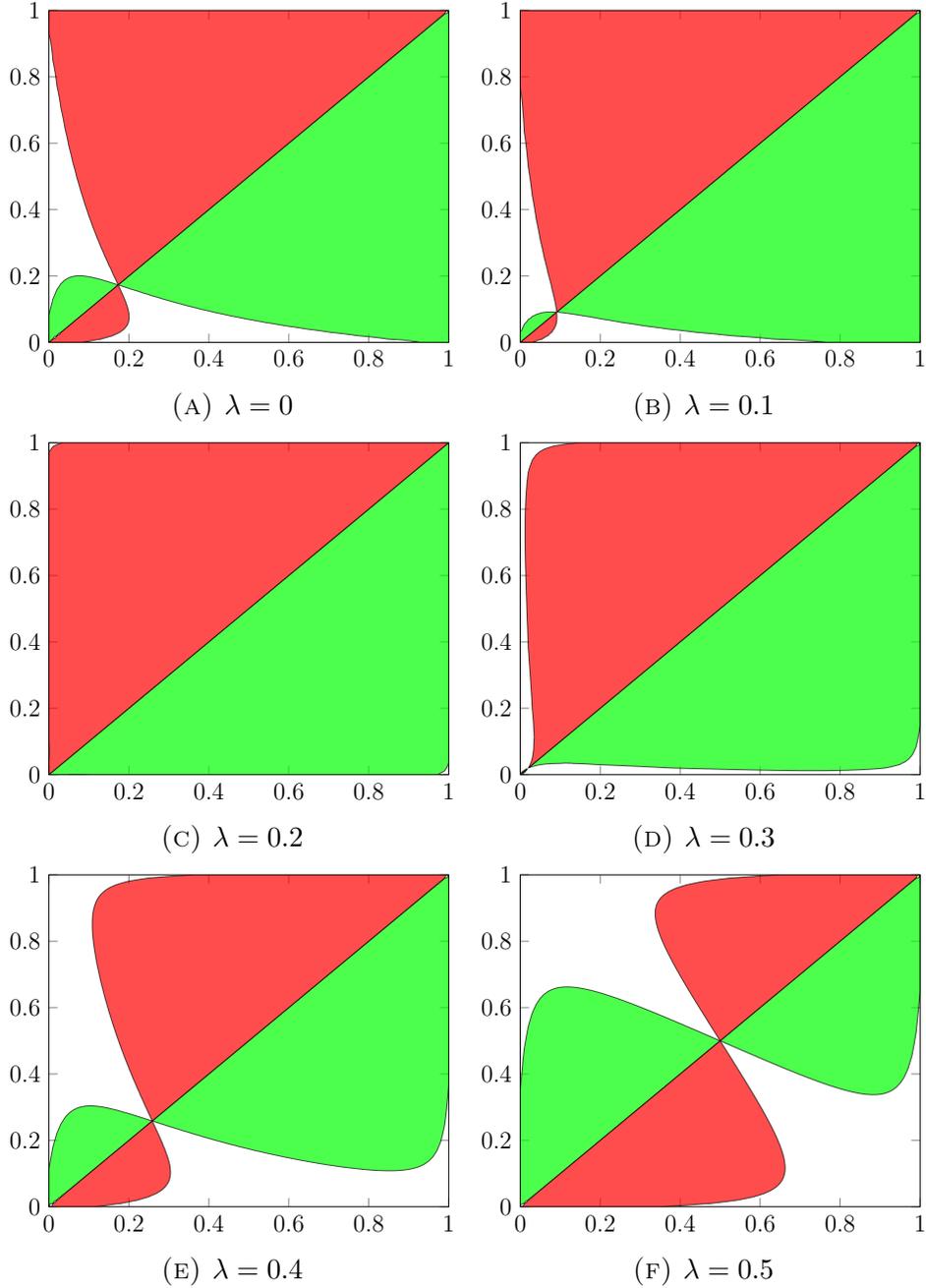

Figure \ref{fig: nodal set double} presents numerical approximations of the nodal sets of $\Lambda(p,q)$ and $\Lambda(q,p)$ for a free growth function of the form $a(x,y)=\lambda A(x)+(1-\lambda)A(y)$, which correspond to the stablity of semi-trival steady states $(N_p, 0)$ and $(0, N_q)$, respectively. The domain $\Omega$ is again a disk of radius $2$ centered at the origin, the anisotropic diffusion parameters $\underline{D}$ and $\overline{D}$ are respectively equal to $0.1$ and $10$, but the function $A$  is now $A(x)=4-\frac{1}{4}\,x^2$. Due to a symmetry of $\Lambda(p,q)$ and $\Lambda(q,p)$ with respect to $\lambda$, we only plot the nodal sets for values of $\lambda$ between $0$ and $0.5$.

For Figure \ref{fig: nodal set double}, in the green colored region $(N_p,0)$ is stable and $(0,N_q)$ is unstable. By Theorems \ref{thm:Fpositive-1-aa}, \ref{thm:Fpositive-2-aa} and \ref{thm:Fchangessign-1}, $(N_p,0)$ is globally stable for $(p,q)$ in the green region. Similarly, the one colored in red is the one in which $(0,N_q)$ is globally stable. The white region is where both $(N_p,0)$ and $(0,N_q)$ are unstable, and there is a unique positive steady state which is globally stable.
The white regions for $\lambda=0.4$ and $\lambda=0.5$ are substantially greater than those for smaller values of $\lambda$. Biologically, this suggests that if the spatial variations of the resource distribution in the vertical and horizontal directions become more comparable, the chances for the coexistence of the two competing populations could be greater. 


\section{The invasion fitness}\label{sec:invasion fitness}

In this section we consider the stability of $(N_p, 0)$,  study some properties of the invasion fitness and establish Theorem~\ref{thm:F-1}.

The linear stability of $(N_p,0)$ is determined by the sign of the smallest eigenvalue, denoted by $\Lambda:=\Lambda(p,q)$, of the linear eigenvalue problem
\begin{equation}\label{eq:varphi-1}
\left\{
\aligned
&D(q)\varphi_{xx}+D(1-q)\varphi_{yy}+(a-N_p)\varphi+\Lambda(p, q)\varphi=0\text{ in }\Omega,\\
&(D(q)\varphi_x, D(1-q)\varphi_y)\cdot \nu=0 \text{ on }\partial\Omega,
\endaligned
\right.
\end{equation}
as seen in the %
following result
.

\begin{lemma}\label{lemma:creteria}
The semi-trivial equilibrium $(N_p,0)$ is linearly stable if $\Lambda(p,q)$ is positive and unstable if $\Lambda(p,q)$ is negative. Similarly, the semi-trivial equilibrium $(0, N_q)$ is linearly stable if $\Lambda(q,p)$ is positive and unstable if $\Lambda(q,p)$ is negative.
\end{lemma}

The proof of Lemma \ref{lemma:creteria} is the same as that of Lemma 5.5 in \cite{CHL2008} and is thus omitted.

\medskip

In the theory of adaptive dynamics \cite{Dieckmann1996, Odo2003, Geritz1998}, $\Lambda(p, q)$ is termed as the \emph{invasion fitness} or \emph{invasion exponent}, which can be regarded as the payoff function for the mutant phenotype with trait $q$, when the resident phenotype with trait $p$ is at the equilibrium. Namely, if $\Lambda(p, q)$ is positive, the mutant with trait $q$ can invade when rare; {on the contrary}, when $\Lambda(p,q)$ is negative, the mutant with trait $q$ cannot invade when rare.
We shall now give some properties of {the map} $(p,q) \mapsto \Lambda(p,q)$ that will be used later to describe more precisely the stability of both semi-trivial steady states.

\subsection{The selection gradient}
If $p=q$, that is when both phenotypes are identical, $\Lambda(p, p)=0$ for any $p$ {in} $[0,1]$, thus both semi-trivial steady states $(N_p,0)$ and $(0, N_q)$ are neutrally stable. In this section, we consider the stability of $(N_p,0)$ for $p$ {and} $q$ sufficiently close to each other. The following result provides a {criterion} in determining the sign of $\Lambda(p,q)$ {in such case} (see also \cite{Slover}).

\begin{lemma} \label{lemma:fitness-formula}
There holds
\[
\left.\frac{\partial\Lambda}{\partial q}\right|_{q=p}=\frac{(\overline{D}-\underline{D}) }{\int_\Omega(N_p)^2}\,F(p),
\]
where $F$ is the function defined by \eqref{eq:functionF}.
\end{lemma}
\begin{proof}
Consider the positive eigenfunction $\varphi$ associated with $\Lambda(p,q)$ and uniquely determined by $\int_\Omega\varphi^2=\int_\Omega(N_p)^2$. It is a smooth function of $p$ and $q$, and, for simplicity of notation, we denote $\frac{\partial\varphi}{\partial q}$ by $\varphi'$. Differentiating system \eqref{eq:varphi-1} with respect to $q$, we obtain
\begin{equation}\label{eq:varphi'}
\left\{
\aligned
D(q)\varphi'_{xx}+D(1-q)\varphi'_{yy}+(a-N_p)\varphi'+\Lambda(p,q)\varphi'+\frac{\partial\Lambda}{\partial q}{(p,q)}\varphi&\\
+(\overline{D}-\underline{D})(\varphi_{xx}-\varphi_{yy})&=0\text{ in }\Omega,\\
(D(q)\varphi'_x+(\overline{D}-\underline{D})\varphi_x,D(1-q)\varphi'_y-(\overline{D}-\underline{D})\varphi_y)\cdot \nu=0\text{ on }\partial\Omega.&
\endaligned
\right.
\end{equation}
Multiplying {the first equation in system \eqref{eq:varphi-1}} by $\varphi'$, integrating {by parts the result over} $\Omega$ and using {the second equation in system \eqref{eq:varphi-1}}, we have
\[
-\int_\Omega(D(q)\varphi_x\varphi'_x+D(1-q)\varphi_y\varphi'_y)+\int_\Omega(a-N_p)\varphi\varphi'+\Lambda{(p,q)}\int_\Omega\varphi\varphi'=0.
\]
Similarly, multiplying {the first equation in system \eqref{eq:varphi'}} by $\varphi$, integrating {by parts the result over} $\Omega$ and using {the second equation in system \eqref{eq:varphi'}}, we obtain
\[
\aligned
-\int_\Omega(D(q)\varphi'_x\varphi_x+D(1-q)\varphi'_y\varphi_y)+\int_\Omega(a-N_p)\varphi'\varphi-(\overline{D}-\underline{D})\int_\Omega((\varphi_x)^2-(\varphi_y)^2)&\\
+\Lambda{(p,q)}\int_\Omega\varphi'\varphi+\frac{\partial \Lambda}{\partial q}{(p,q)}\int_\Omega\varphi^2&=0.
\endaligned
\]
Subtracting the above two {equalities then yields}
\begin{equation}\label{eq:formula-2-1}
\frac{\partial \Lambda}{\partial q}{(p,q)}\int_\Omega\varphi^2=(\overline{D}-\underline{D})\int_\Omega((\varphi_x)^2-(\varphi_y)^2).
\end{equation}
{Since} $\Lambda(p,p)=0$, {it follows from the normalization of $\varphi$ that} $\varphi_{|_{q=p}}=N_p$, which completes the proof.
\end{proof}

In view of Lemma \ref{lemma:fitness-formula}, it is critical to understand the sign of function $F$. The following result establishes Theorem \ref{thm:F-1}.

\begin{lemma}\label{lem:selection-gradient-1}
For any $\theta$ {in} $[0,1]$, {one has} $F'(\theta)<0$.
\end{lemma}
\begin{proof}
For {simplicity of notation, } denote $\frac{\partial N_\theta}{\partial\theta}$ by $N'_{\theta}$. {Integrating by parts, one gets}
\[
\aligned
F'(\theta)&=2\int_\Omega ( (N_\theta)_x (N'_\theta)_x-(N_\theta)_y (N'_\theta)_y )\\
&=2\int_{\partial\Omega} N'_{\theta} ((N_\theta)_x, -(N_\theta)_y)\cdot \nu-2\int_{\Omega} N'_{\theta} ((N_\theta)_{xx}-(N_\theta)_{yy}).
\endaligned
\]
Differentiating system in \eqref{eq:singlespeciessteady} for $N_\theta$ with respect to $\theta$, we obtain
\[
\left\{
\aligned
&D(\theta)(N'_\theta)_{xx}+D(1-\theta)(N'_\theta)_{yy}+(a-2N_\theta)N'_{\theta}+(\overline{D}-\underline{D})((N_\theta)_{xx}-(N_\theta)_{yy})=0\text{ in }\Omega,\\
&(D(\theta)(N'_\theta)_{x}+(\overline{D}-\underline{D})(N_\theta)_{x},D(1-\theta)(N'_\theta)_{y}-(\overline{D}-\underline{D})(N_\theta)_{y})\cdot\nu=0\text{ on }\partial\Omega.
\endaligned
\right.
\]
Multiplying the {first of the above equations} by $N'_\theta$ and integrating the result {over} $\Omega$, we {find that}
\[
\aligned
&(\overline{D}-\underline{D})\int_\Omega N'_\theta((N_\theta)_{xx}-(N_\theta)_{yy})\\
&=-\int_\Omega N'_\theta\left[D(\theta)(N'_\theta)_{xx}+D(1-\theta)(N'_\theta)_{yy}
+(a-2N_\theta)N'_{\theta}\right]\\
&=-\int_{\partial\Omega} N'_\theta (D(\theta)(N'_\theta)_{x},
D(1-\theta) (N'_\theta)_{y})\cdot\nu\\
&\quad+\int_\Omega [D(\theta)((N'_\theta)_{x})^2+D(1-\theta)((N'_\theta)_{y})^2- (a-2N_\theta)(N'_\theta)^2]\\
&=(\overline{D}-\underline{D})\int_{\partial\Omega} N'_\theta ((N_\theta)_{x},-(N_\theta)_{y})\cdot \nu\\
&\quad+\int_\Omega[D(\theta)((N'_\theta)_{x})^2+D(1-\theta)((N'_\theta)_{y})^2-(a-2N_\theta)(N'_\theta)^2],
\endaligned
\]
where the last equality follows from the boundary condition {satisfied by} $N'_\theta$. {We therefore have}
\[
F'(\theta)=-\frac{2}{\overline{D}-\underline{D}}\int_\Omega[D(\theta)((N'_\theta)_{x})^2+D(1-\theta)((N'_\theta)_{y})^2-(a-2N_\theta)(N'_\theta)^2].
\]

Let $\lambda_1$ denote the smallest eigenvalue of the linear problem
\[
\left\{
\aligned
&D(\theta)\varphi_{xx}+D(1-\theta)\varphi_{yy}+(a-N_\theta)\varphi+\lambda\varphi=0\text{ in }\Omega,\\
&(D(\theta)  \varphi_x , D(1-\theta) \varphi_y)\cdot\nu=0\text{ on }\partial\Omega.
\endaligned
\right.
\]
It is well-known that $\lambda_1$ can be characterized by the variational formula
\begin{equation}\label{eq:temp1}
\lambda_1=\inf_{\varphi\not=0,\ \varphi\in H^1(\Omega)}\frac{\int_\Omega[D(\theta)(\varphi_x)^2+D(1-\theta)(\varphi_y)^2-(a-N_\theta) \varphi^2]}{\int_\Omega \varphi^2},
\end{equation}
and that  $\lambda_1$ is the only eigenvalue such that its corresponding eigenfunction does not change sign in $\Omega$. {Using system \eqref{eq:singlespeciessteady} for $N_\theta$, we infer} that $\lambda_1=0$ and {that} its corresponding eigenfunction is a scalar multiple of $N_\theta$. In particular, by choosing the test function $\varphi=N'_\theta$ in \eqref{eq:temp1}, we {have} that
\[
\int_\Omega [D(\theta)((N'_\theta)_x)^2+D(1-\theta)((N'_\theta)_y)^2-(a-N_\theta)(N'_\theta)^2]\ge 0,
\]
{which yields}
\[
F'(\theta)\le-\frac{2}{\overline{D}-\underline{D}}\int_\Omega N_\theta(N'_\theta)^2\le 0.
\]
{This gives} $F'\le 0$, {with the equality} if and only if $N'_\theta\equiv 0$ in $\Omega$.

Finally, 
if $N'_\theta\equiv 0$ for some $\theta$, then $N_\theta$ satisfies
\[
\left\{
\aligned
& (N_\theta)_{xx}-(N_\theta)_{yy}=0\text{ in }\Omega,\\
&((N_\theta)_{x},-(N_\theta)_{y})\cdot \nu=0\text{ on }\partial\Omega.
\endaligned
\right.
\]
This, together with the boundary condition {satisfied by} $N_\theta$, implies that
\[
(N_\theta)_{x}\nu_x=(N_\theta)_{y}\nu_y=0\text{ on }\partial\Omega.
\]
{According to} Lemma \ref{lemma:wave}
in the Appendix, $N_\theta$ is {then} a positive constant function, and {it follows from the first equation in \eqref{eq:singlespeciessteady} that the function} $a(x,y)$ is also a constant function, which contradicts assumption {(A1)}. {Consequently, one has} $F'<0$ in $[0, 1]$.
\end{proof}

\begin{remark}
If we allow Lipschitz domains with flat parts on the boundary, it is possible to construct domains such that $F\equiv 0$ (and thus $F'\equiv 0$) in $[0, 1]$. See {the} Appendix \ref{appendix:2} for further discussions.
\end{remark}







\subsection{Concavity of $\Lambda(p,q)$}
The following result concerns the concavity of $\Lambda(p, q)$ with respect to $q$.

\begin{lemma}\label{lemma:concavity}
{For any $p$ in $[0,1]$}, the function $q\mapsto\Lambda(p,q)$ is concave on $[0,1]$. Moreover, if $\Lambda_{q}(p^*,q^*)=0$ for some {couple} $(p^*,q^*)$, then $\Lambda_{qq}(p^*,q^*)<0$.
\end{lemma}
\begin{proof}
{Let us} fix $p$ {in} $[0,1]$. The concavity of $q\mapsto\Lambda(p,q)$ follows from {a} standard argument based on {the} variational characterization of $\Lambda(p,q)$, see \cite{Ni2011}. {Nevertheless}, we include {here a proof of this result in order} to facilitate the proof of {the second statement of the Lemma}.

{Consider the positive eigenfunction $\varphi$ associated with $\Lambda(p,q)$  such that $\int_\Omega\varphi^2=\int_\Omega(N_p)^2$.}
Differentiating \eqref{eq:formula-2-1} with respect to $q$ yields
\[
\Lambda_{qq}{(p,q)}\int_\Omega \varphi^2+2\Lambda_q{(p,q)}\int_\Omega\varphi'\varphi=2(\overline{D}-\underline{D})\int_\Omega(\varphi'_x\varphi_x-\varphi'_y\varphi_y),
\]
{with the notations $\varphi'=\frac{\partial\varphi}{\partial q}$,} $\varphi'_x=\frac{\partial\varphi_x}{\partial q}$, and $\varphi'_y=\frac{\partial\varphi_y}{\partial q}$. {Note that the second term in the left-hand side vanishes due to the normalization condition on $\varphi$, which implies that  $\int_\Omega\varphi'\varphi=0$.}

Multiplying {the first equation in system} \eqref{eq:varphi'} by $\varphi'$ and integrating {by parts} the result {over} $\Omega$ {then gives}
\[
\aligned
&(\overline{D}-\underline{D})\int_\Omega (\varphi_x \varphi'_x-\varphi_y \varphi'_y)\\
&=-\int_\Omega\left[D(q)(\varphi'_x)^2+D(1-q)(\varphi'_y)^2 -(a-N_p)(\varphi')^2
-\Lambda{(p,q)}(\varphi')^2\right]\\
&\le 0,
\endaligned
\]
where we {have used the second equation in system \eqref{eq:varphi'} and the fact that} $\int_\Omega\varphi'\varphi=0$, the inequality {following} from the variational characterization of $\Lambda(p,q)$,
\[
\Lambda(p,q)=\inf_{\varphi\not=0,\ \varphi\in H^1(\Omega)}
\frac{\int_\Omega\left[D(q)(\varphi_x)^2+D(1-q)(\varphi_y)^2-(a-N_p)\varphi^2\right]}{\int_\Omega \varphi^2}.
\]
{It then holds that} $\Lambda_{qq}{(p,q)}\le 0$,
{where the equality holds} if and only if the function $\varphi'$ is a scalar multiple of $\varphi$. Since $\varphi$ is positive in $\Omega$, {this implies that $\varphi'\equiv 0$  in $\Omega$.}

Thus, if $\Lambda_q{(p,q)}=\Lambda_{qq}{(p,q)}=0$, {system \eqref{eq:varphi'} reduces to}
\[ 
\varphi_{xx}-\varphi_{yy}=0\text{ in }\Omega,\quad \text{ and }\quad (\varphi_x,-\varphi_y)\cdot\nu=0\text{ on }\partial\Omega.
\] 
Due to the boundary condition for $\varphi$ in system \eqref{eq:varphi-1}, we further have $\varphi_x\nu_x=\varphi_y \nu_y=0$ on $\partial\Omega$. 
As a consequence of Lemma \ref{lemma:wave}, the function $\varphi$ is constant, which implies, {using the first equation in system \eqref{eq:varphi-1}, that the function} $a-N_p$ is also constant. Integrating {over $\Omega$} the main equation in system \eqref{eq:singlespeciessteady} with $\theta=p$, we then obtain
\[
\int_\Omega(a-N_p)N_p=0,
\]
which ensures that $a-N_p=0$ in $\Omega$, so that system \eqref{eq:singlespeciessteady} for $N_p$ reduces to
\[ 
\left\{
\aligned
&D(p)(N_p)_{xx}+D(1-p)(N_p)_{yy}=0\quad \ \ \text{ in }\Omega,\\
&(D(p)(N_p)_{x},D(1-p)(N_p)_{y})\cdot\nu=0\quad \text{ on }\partial\Omega.
\endaligned
\right.
\] 
{It finally follows from} the maximum principle {that the function} $N_p$ is constant. This contradicts the assumption {of the function} $a$ being non-constant.
\end{proof}

\begin{lemma}\label{lemma:nondegeneracy}
If $\Lambda(\tilde{p},\tilde{q})=0$ for some {couple} $(\tilde{p},\tilde{q})$ {such that} $\tilde{p}\not=\tilde{q}$, then $\Lambda_q(\tilde{p},\tilde{q})\not=0$.
\end{lemma}
\begin{proof}
We argue by contradiction. Suppose that there exists {a couple $(\tilde{p},\tilde{q})$} {such that} $\tilde{p}\not=\tilde{q}$ {for which} $\Lambda(\tilde{p},\tilde{q})=\Lambda_q(\tilde{p},\tilde{q})=0$. By Lemma \ref{lemma:concavity}, {one has} $\Lambda_{qq}(\tilde{p},\tilde{q})<0$, which implies that {the function $q\mapsto\Lambda(\tilde{p},q)$} {has a} local maximum {point} at $q=\tilde{q}$. This contradicts the fact that $\Lambda(\tilde{p},\tilde{p})=0$ and the concavity of $\Lambda{(p,q)}$ in $q$.
\end{proof}

A consequence of Lemma \ref{lemma:nondegeneracy} is that the nodal set of $\Lambda{(p,q)}$ within the region $\{(p,q)\ :\ 0\le p,q\le1\}$ can be parameterized by a function $q=\tilde{q}(p)$. However, we caution the reader that the domain of {this} function can be either $[0,1]$ or a strict subset of {it}.

\subsection{Numerical results}
Numerical approximations of the graph of the function $F$ in the case of a free-growth function of the form $a(x,y)=\lambda A(x)+(1-\lambda)A(y)$ are provided in Figure \ref{fig:F}. For the simulations, we considered a disk of radius $2$ centered at the origin for the domain $\Omega$, anisotropic diffusion parameters $\underline{D}$ and $\overline{D}$ respectively equal to $0.1$ and $10$, and the function $A(x)=2-\sin(\pi x)$. Due to a symmetry\footnote{It is easily seen that, for all $\lambda$ in $[0,1]$ and all $\theta$ in $[0,1]$, the value $F(\theta)$ for $\lambda$ is equal to the value of $-F(1-\theta)$ for $1-\lambda$.} in the function $F$ with respect to $\lambda$, we only plot graphs for different values of $\lambda$ between $0$ and $0.5$, illustrating how the function $F$ goes from strictly negative, to sign-changing once, and to strictly positive as $\lambda$ varies.

\begin{figure}
\centering
\subfigure[$\lambda=0$]{\resizebox{0.45\textwidth}{!}{\begin{tikzpicture}
\begin{axis}[xmin=0,xmax=1]
\addplot[] table {case1_l0F.txt} ;
\end{axis}
\end{tikzpicture}}}
\subfigure[$\lambda=0.1$]{\resizebox{0.45\textwidth}{!}{\begin{tikzpicture}
\begin{axis}[xmin=0,xmax=1]
\addplot[] table {case1_l2F.txt} ;
\end{axis}
\end{tikzpicture}}}

\subfigure[$\lambda=0.2$]{\resizebox{0.45\textwidth}{!}{\begin{tikzpicture}
\begin{axis}[xmin=0,xmax=1]
\addplot[] table {case1_l4F.txt} ;
\end{axis}
\end{tikzpicture}}}
\subfigure[$\lambda=0.3$]{\resizebox{0.45\textwidth}{!}{\begin{tikzpicture}
\begin{axis}[xmin=0,xmax=1]
\addplot[] table {case1_l6F.txt} ;
\end{axis}
\end{tikzpicture}}}

\subfigure[$\lambda=0.4$]{\resizebox{0.45\textwidth}{!}{\begin{tikzpicture}
\begin{axis}[xmin=0,xmax=1]
\addplot[] table {case1_l8F.txt} ;
\end{axis}
\end{tikzpicture}}}
\subfigure[$\lambda=0.5$]{\resizebox{0.45\textwidth}{!}{\begin{tikzpicture}
\begin{axis}[xmin=0,xmax=1]
\addplot[] table {case1_l10F.txt} ;
\end{axis}
\end{tikzpicture}}}
\caption{Numerical approximations of the graphs of the function $F$ for $\lambda$ taking the values $0$, $0.1$, $0.2$, $0.3$, $0.4$, and $0.5$, in the case where $\Omega$ is a disk of radius $2$ centered at the origin, $\underline{D}=0.1$, $\overline{D}=10$, and $A(x)=2-\sin(\pi x)$.}\label{fig:F}
\end{figure}
\emd{comment}

\section{Local stability for $(p,q)$ in $[0,1]^2$: nodal set of invasion fitness}\label{sec:local stability}
{We have previously} considered the stability of {the semi-trivial steady state} $(N_p,0)$ for $p$ {and} $q$ close to each other. In this section, we {study} the local stability of $(N_p,0)$ for general $p$ {and} $q$ in $[0,1]$. This is equivalent to giving a description of the nodal set of $\Lambda(p,q)$ for $p$ {and} $q$ {in} $[0,1]$, and as well the sets where $\Lambda(p,q)$ is positive or negative.

\smallskip

By Lemma \ref{lem:selection-gradient-1}, {the function} $F$ is strictly decreasing in $[0, 1]$, so it suffices to consider three generic cases: $F(\theta)>0$ for all $\theta$ {in} $[0,1]$, $F$ changes sign exactly once in $(0,1)$, and $F(\theta)<0$ for all $\theta$ {in} $[0,1]$.

\subsection{Stability when $F(\theta)>0$}
First, we consider the case {for which the function} $F$ {is positive} in $[0,1)$.

\begin{lemma}\label{lemma:4.1}
Suppose that $F(\theta)>0$ for all $\theta$ {in} $[0,1)$. If $0\le p<q\le 1$, then {the steady state} $(0, N_q)$ is unstable. {Conversely}, if $0\le q<p\le 1$, then {the steady state} $(N_p, 0)$ is unstable.
\end{lemma}
\begin{proof}
We first consider the situation $0\le p<q<1$. The stability of {the steady state} $(0, N_q)$ is determined by the sign of the smallest eigenvalue, denoted by $\Lambda(q,p)$, of the linear eigenproblem
\[
\left\{
\aligned
&D(p)\varphi_{xx}+D(1-p)\varphi_{yy}+(a-N_q)\varphi+\lambda\varphi=0\text{ in }\Omega,\\
&(D(p)\varphi_x, D(1-p)\varphi_y)\cdot{\nu}=0\text{ on }\partial\Omega.
\endaligned
\right.
\]
By the variational characterization of $\Lambda(q, p)$ we have
\[
\aligned
\Lambda(q,p)
&=\inf_{\varphi\not=0,\ \varphi\in H^1(\Omega)}
\frac{\int_\Omega \left[D(p)\varphi_x^2+D(1-p)\varphi_y^2-(a-N_q)\varphi^2\right]}{\int_\Omega\varphi^2}\\
&\le\frac{\int_\Omega[D(p)((N_q)_x)^2+D(1-p)((N_q)_y)^2-(a-N_q)(N_q)^2]}{\int_\Omega(N_q)^2}
\endaligned
\]

Multiplying the main equation in system \eqref{eq:singlespeciessteady} with $\theta=q$ by $N_q$ and integrating {by parts} the result {over} $\Omega$ {yields}
\[
\int_\Omega\left[D(q)\left((N_q)_x\right)^2+D(1-q)\left((N_q)_y\right)^2-(a-N_q)(N_q)^2\right]=0.
\]

Therefore, we have, {using the assumptions that $p<q$ and $F>0$ in $[0,1]$},
\[
\aligned
\Lambda(q,p)
&\le \frac{[D(p)-D(q)]\int_\Omega((N_q)_x)^2+[D(1-p)-D(1-q)]\int_\Omega((N_q)_y)^2}{\int_\Omega(N_q)^2}\\
&=(\overline{D}-\underline{D})(p-q)\frac{\int_\Omega ((N_q)_x)^2-\int_\Omega ((N_q)_y)^2}{\int_\Omega (N_q)^2}\\
&=(\overline{D}-\underline{D})(p-q)\frac{F(q)}{\int_\Omega(N_q)^2}<0.
\endaligned
\]
{If} $q=1$, we note that, as $\Lambda(p,q)<0$ for $0\le p<q<1$, $\Lambda(p,1)\le 0$. Since $\Lambda(1,1)=0$, we see by Lemma \ref{lemma:nondegeneracy} that $\Lambda(p,1)<0$ for all $0\le p<1$. Hence, {the steady state $(0,N_q)$} is unstable for $0\le q<p\le 1$.

Similarly, we can show that if $0\le q<p\le 1$, then 
$(N_p,0)$ is unstable.
\end{proof}

The first main result of this section is {the following.}

\begin{theorem}\label{thm:Fpositive-1}
Suppose {that the function $F$ is positive} in $[0,1]$. Then, there exists some positive, continuous function $q=q^*(p)$, defined in $[0,1]$, satisfying  $p<q^*(p)\le 1$ for {all $p$ in} $[0,1]$ and $q^*(p)\equiv 1$ for $p$ close to $1$, such that
\begin{enumerate}[\rm(i)]
\item $\Lambda(p,q)>0$ for $0\le p<q<q^*(p)$,\label{thm:Fpositive-1 - statement positive}
\item $\Lambda(p,q)<0$ for
$q^*(p)<q\le 1$ and $0\le q<p\le 1$.\label{thm:Fpositive-1 - statement negative}
\end{enumerate}
In particular, if $q^*(p)\equiv 1$, then $\Lambda(p,q)>0$ for $p<q\le 1$ and $\Lambda(p,q)<0$ for $0\le q<p$.
\end{theorem}
\begin{proof}
Since $\Lambda(p, p)=0$ for all $p$ {in} $[0,1]$ and $F$ {is positive} in $[0, 1]$, there exists some {positive real number} $\delta$ such that $\Lambda(p,q)>0$ for {$p$ and $q$ in $[0,1]$ with $0<q-p<\delta$}.

Given any $p$ {in} $[0,1]$, if $\Lambda(p,1)<0$, by Lemma \ref{lemma:concavity} and the positivity of $\Lambda(p, q)$ in the strip $0<q-p<\delta$, there exists a unique $q^*=q^*(p)$ {in} $(p, 1)$ such that $\Lambda(p, q^*)=0$, $\Lambda(p, q)>0$ for $p<q<q^*$ and $\Lambda(p, q)<0$ for $q>q^*$. If $\Lambda(p, 1)\ge 0$, we define $q^*(p)=1$. Again by Lemma \ref{lemma:concavity}, $\Lambda(p, q)>0$ holds for $p<q<q^*=1$. This establishes {statements} \eqref{thm:Fpositive-1 - statement positive} and \eqref{thm:Fpositive-1 - statement negative}.

To show that $q=q^*(p)$ is a continuous curve, suppose that $\Lambda(\tilde{p}, \tilde{q})=0$ for some $\tilde{p}<\tilde{q}<1$. By Lemma \ref{lemma:nondegeneracy}, there exists a smooth curve $q=q^{**}(p)$ passing through $(\tilde{p},\tilde{q})$ {and} such that $\Lambda(p,q^{**}(p))=0$, {which} can be extended to the left and {to} right until it reaches either $p=0$ or $q=1$, as by the choice of $\delta$, $q=q^{**}(p)$ and its extension can never enter the strip $0<q-p<\delta$. For each $p$, there exists at most one $q>p$ such that $\Lambda(p, q)=0$.

Hence, $q^{**}\equiv q^*$ as long as these functions are strictly less than $1$. Therefore, $q^*$ defines a continuous curve on $[0,1]$. By the choice of $\delta$, we see that $q^*(p)\equiv 1$ for $p$ close to {$1$}.
\end{proof}

Similarly to Theorem \ref{thm:Fpositive-1}, the local stability of {the semi-trivial steady state} $(0,N_q)$ can be determined as follows.

\begin{theorem}\label{thm:Fpositive-1-a}
Suppose {that the function $F$ is positive} in $[0,1]$. Then, there exists some positive, continuous function $p=p^*(q)$, defined in $[0,1]$, satisfying  $q<p^*(q)\le 1$ for $q$ in $[0,1]$ and $p^*(q)\equiv 1$ for $q$ close to $1$, such that $\Lambda(q,p)>0$ for $q<p<p^*(q)$ and $\Lambda(q,p)<0$ for $p^*(q)<p\le 1$ and $q<p$.
\end{theorem}

The proof being the same as {that} for Theorem \ref{thm:Fpositive-1}, we omit it.

\subsection{Stability when $F(\theta)<0$}
Next, we consider the case {for which the function} $F$ {is negative} in $[0,1]$. This case is similar to the previous one, so that we may state the following results without proofs.

\begin{theorem}\label{thm:Fpositive-2}
Suppose {that the function $F$ is negative} in $[0, 1]$. Then, there exists {a} positive continuous function $q=q^*(p)$, defined in $[0, 1]$, satisfying  $0\le q^*(p)<p$ for $p$ in $[0,1]$ and $q^*(p)\equiv 0$ for $p$ close to $0$ such that
\begin{enumerate}[\rm (i)]
\item $\Lambda(p,q)>0$ for $q^*(p)<q<p\le 1$,
\item $\Lambda(p,q)<0$ for 
$0\le q<q^*(p)$ and  $0\le p<q\le 1$.
\end{enumerate}
In particular, if $q^*(p)\equiv 0$, then $\Lambda(p,q)>0$ for $q<p$ and $\Lambda(p,q)<0$ for $q>p$.
\end{theorem}


\begin{theorem}\label{thm:Fpositive-2-a}
Suppose that the function $F$ is negative in $[0,1]$. Then, there exists {a} positive continuous function $p=p^*(q)$, defined in $[0,1]$, satisfying  $p^*(q)<q$ for $q$ in $[0,1]$ and $p^*(q)\equiv 0$ for $q$ close to $0$, such that $\Lambda(q,p)>0$ for $p^*(q)<p<q$ and $\Lambda(q,p)<0$ for $0\le p<p^*(q)$ and $q<p$.
\end{theorem}


\subsection{Stability when $F$ changes sign once}
We finally consider the case {for which the function} $F$ {possesses} a unique root in $(0,1)$, denoted by $\theta^*$. {This function being} decreasing, {this implies} that {it is positive} in $[0,\theta^*)$ and {negative} in $(\theta^*,1]$.

\begin{theorem}\label{thm:Fpositive-3}
Suppose {that the function} $F$ has a unique root $\theta^*$ in $(0,1)$. Then, there exists {a} positive continuous function $q=q^*(p)$, defined in $[0, 1]$, satisfying  $p<q^*(p)\le 1$ for $p$ {in} $[0,\theta^*)$ and $0\le q^*(p)<p$ for $p$ {in} $(\theta^*,1]$, such that

\begin{enumerate}[\rm(i)]
\item $\Lambda(p,q)>0$ for $p$ {in} $[0,1]$ and $\min\{q^*(p), p\}<q<\max\{q^*(p), p\}$,\label{thm:Fpositive-3 - statement positive}
\item $\Lambda(p,q)<0$ for $p$ {in} $[0,1]$, $\max\{q^*(p), p\}<q\le 1$ and $0\le q<\min\{q^*(p),p\}$.\label{thm:Fpositive-3 - statement negative}
\end{enumerate}
\end{theorem}
\begin{proof}
Since $\Lambda(p, p)=0$ for all $p$ {in} $[0,1]$ and $F$ changes sign exactly once at some $\theta^*$ {in} $(0,1)$, {it follows from} Lemma \ref{lemma:concavity} and the implicit function theorem {that} there exist {both} a smooth curve $q=q^*(p)$, which passes through $(\theta^*, \theta^*)$, such that $q^*(p)>p$ for $p$ less than and close to $\theta^*$, and $q^*(p)<p$ for $p$ greater than and close to $\theta^*$, {and a positive real number} $\delta$, {such that} $\Lambda(p,q)=0$ in the stripe $|q-p|<\delta$ if and only if either $q=p$ or $q=q^*(p)$.

{Using} Lemma \ref{lemma:nondegeneracy}, we can extend the curve $q=q^*(p)$ to the left until it reaches either $p=0$ or $q=1$. Note that we can choose $\delta$ so small that this extension can never {re-enter} the strip $0<q-p<\delta$ once {it leaves it}. If {it} first reaches $q=1$ at some $p=\hat{p}$ {in} $(0,\theta^*)$, one can proceed as in the proof of Theorem \ref{thm:Fpositive-1} to define $q^*(p)$ for $p$ {in} $[0,\theta^*)$ and show that $q=q^*(p)$ is a continuous curve. If {it} never {attains} $q=1$, note that {it} can only intersect the line $q=p$ at $q=p=\theta^*$, thus it can be defined at $p=0$ in such a way that $p<q^*(p)<1$ for $p$ {in} $(0,\theta^*)$. {From} Lemma \ref{lemma:nondegeneracy} and the implicit function theorem, {the curve} $q=q^*(p)$ is smooth in this scenario. Similarly, one can extend $q^*(p)$ to $[\theta^*, 1]$ as a continuous curve.

{These arguments also show} that the nodal set of $\Lambda(p, q)$ is contained in {the line} $q=p$ and {the curve} $q=q^*(p)$, so that the conclusions in \eqref{thm:Fpositive-3 - statement positive} and \eqref{thm:Fpositive-3 - statement negative}
hold.
\end{proof}


Similarly to Theorem \ref{thm:Fpositive-3}, the local stability of  $(0, N_q)$ can be determined as follows.

\begin{theorem}\label{thm:Fpositive-4}
Suppose {that the equation} $F=0$ has a unique root, denoted by $\theta^*$, in $(0, 1)$. Then, there exists some positive, continuous function $p=p^*(q)$, defined in $[0, 1]$, satisfying  $q<p^*(q)\le 1$ for $q\in [0, \theta^*)$ and $0\le p^*(q)<q$ for $q\in (\theta^*, 1]$, such that
\begin{enumerate}[\rm(i)]
\item $\Lambda(q, p)>0$ for $q\in [0, 1]$ and $\min\{p^*(q), q\}<p<\max\{p^*(q), q\}$,
\item $\Lambda(q, p)<0$ for $q\in [0, 1]$, $\max\{p^*(q), q\}<p\le 1$ and $0\le p<\min\{p^*(q), q\}$.
\end{enumerate}
\end{theorem}

The proof being the same as {that} for Theorem \ref{thm:Fpositive-3}, we omit it.

\subsection{Proofs of Theorem \ref{thm:Fpositive-001} and Corollary \ref{thm:F}}{We are now in a position to prove some main results of the paper.}
Theorem \ref{thm:Fpositive-001} follows from Theorems \ref{thm:Fpositive-1}, \ref{thm:Fpositive-2}, and~\ref{thm:Fpositive-3}.

\begin{proof}[Proof of Corollary {\rm\ref{thm:F}}]
If {the function} $F$ {is positive} in $[0,1]$, {statement} \eqref{thm:F - statement positive} follows from Theorem \ref{thm:Fpositive-1}. If $F>0$ in $[0,1)$ with $F(1)=0$, we can apply Lemma \ref{lemma:4.1} to conclude that $p=1$ is an evolutionarily singular strategy but not {an} evolutionarily stable {one}, and $p=0$ is thus the only ESS. {Statement} \eqref{thm:F - statement root} can be proved similarly. Finally, {statement} \eqref{thm:F - statement negative} is a direct consequence of Theorem \ref{thm:Fpositive-3}.
\end{proof}

\begin{remark}\label{remark:11}
While the function $F$ plays a critical role in the analysis provided in the current section, it appears that it only captures some partial information on $\Lambda(p,q)$ and cannot possibly determine {entirely} the nodal set of $\Lambda(p,q)$. For instance, even the sign of $\Lambda(0,1)$ cannot be {resolved using only the} function $F$, as it depends on $a$, $\underline{D}$ and $\overline{D}$ in delicate manners. As an example, assume that $a(x,y)=\lambda A(x)+(1-\lambda)A(y)$, where $\lambda$ belongs to $[0,1)$ and $A$ attains a strict global maximum. Then, for large $\overline{D}$, choosing $\underline{D}$ sufficiently small, $\Lambda(0,1)<0$. However, for such choices of $a$, $\underline{D}$ and $\overline{D}$, the function $F$ changes from negative to sign-changing and to positive as $\lambda$ varies from $0$ to $1$. We refer to the Appendix \ref{appendix:3}
for further  discussions.
\end{remark}

\section{Full dynamics of the two-species model}\label{sec:full dynamics}
\subsection{Local stability of semi-trivial steady states}
In this subsection, we {investigate further} the local stability of both semi-trivial steady states $(N_p,0)$ and $(0,N_q)$, for general $p$ and $q$ in $[0,1]$.

The {next} result {shows} that {the states} $(N_p,0)$ and $(0,N_q)$ cannot be {simultaneously stable}, \textit{i.e.} bistability cannot occur.

\begin{lemma}\label{lem:bistability-1}
{The following assertions hold for any $p$ and $q$ in $[0,1]$.}
\begin{enumerate}[(i)]
\item If $\Lambda(p,q)>0$, then $\Lambda(q,p)<0$.\label{lem:bistability-1 - statement 1}
\item If $\Lambda(p,q)=0$, then either $p=q$ or $\Lambda(q, p)<0$.\label{lem:bistability-1 - statement 2}
\end{enumerate}
\end{lemma}
\begin{proof}
{It follows from} the variational characterization of $\Lambda(p,q)$ {that}
\[
\aligned
\Lambda(p,q)
&=\inf_{\varphi\not=0,\ \varphi\in H^1(\Omega)}\frac{\int_\Omega\left[D(q)(\varphi_x)^2+D(1-q)(\varphi_y)^2-(a-N_p)\varphi^2\right]}{\int_\Omega \varphi^2}\\
&\le\frac{\int_\Omega[D(q)((N_q)_x)^2+D(1-q)((N_q)_y)^2-(a-N_p)(N_q)^2]}{\int_\Omega(N_p)^2}\\
&=\frac{\int_\Omega (N_p-N_q)(N_q)^2}{\int_\Omega(N_p)^2},
\endaligned
\]
where the last equality follows from the equation of $N_q$. If $\Lambda(p,q)>0$,
then
\[
\int_\Omega(N_q)^3<\int_\Omega N_p(N_q)^2\le\left(\int_\Omega(N_p)^3\right)^{1/3}\left(\int_\Omega(N_q)^3\right)^{2/3},
\]
which implies that
\[
\int_\Omega (N_q)^3<\int_\Omega (N_p)^3.
\]
If we assume that $\Lambda(q,p)\ge 0$, we have, by the same argument as above,
\[
\int_\Omega (N_p)^3\le \int_\Omega (N_p)^2 N_q
\le \left(\int_\Omega (N_p)^3\right)^{2/3}\left(\int_\Omega (N_q)^3\right)^{1/3},
\]
from which we get
\[
\int_\Omega (N_p)^3\le \int_\Omega (N_q)^3,
\]
which is a contradiction. Hence, {statement} \eqref{lem:bistability-1 - statement 1} holds.

If $\Lambda(p, q)=0$, following the same argument as above, we see that $\Lambda(q, p)\le 0$. If $\Lambda(q, p)=0$, the only possibility is that $N_p\equiv N_q$, {that is} $p=q$, which proves 
\eqref{lem:bistability-1 - statement 2}.
\end{proof}

To describe the global dynamics of the two-species model, we first introduce the sets 
\[
\aligned
\Gamma_1:&=\{(p,q)\in [0,1]^2: \Lambda(p, q)=0\},\\ 
\Gamma_2:&=\{(p,q)\in [0,1]^2: \Lambda(q, p)=0\}. 
\endaligned
\]
Clearly, $\Gamma_1$ and $\Gamma_2$ correspond to 
the situation when $(N_p, 0)$ and $(0, N_q)$ are neutrally stable. Therefore, we have, by Theorems \ref{thm:Fpositive-1} to~\ref{thm:Fpositive-4},
{
\[
\aligned
&\Gamma_1=\{(p,q)\in [0,1]^2: q=p\text{ or }q=q^*(p)\} 
, 
\\
&\Gamma_2=\{(p,q)\in [0,1]^2: p=q\text{ or }p=p^*(q)\} 
.
\endaligned
\]
}

Next, we define the sets 
{
\[
\aligned
\Sigma_1:&=\{(p,q)\in [0,1]^2: \Lambda(p,q)>0>\Lambda(q,p)\},\\
\Sigma_2:&=\{(p,q)\in [0,1]^2: \Lambda(p,q)<0<\Lambda(q,p)\},\\
\Sigma_3:&=\{(p,q)\in [0,1]^2: \Lambda(p,q)<0,\ \Lambda(q,p)<0\}.
\endaligned
\]
The sets} $\Sigma_i$ ($i=1,2,3$) are disjoint and
\[
\Sigma_1\cup \Sigma_2\cup\Sigma_3=[0,1]\times[0,1]/(\Gamma_1\cup\Gamma_2\cup\Gamma_3).
\]

\begin{theorem}\label{thm:bistability-2}
The following characterizations hold:
\[
\aligned
\Sigma_1&=\{(p,q)\in [0,1]\times [0,1]: (q-q^*(p))(q-p)<0\},\\
\Sigma_2&=\{(p,q)\in [0,1]\times [0,1]: (p-p^*(q))(p-q)<0\},\\
\Sigma_3&=\{(p,q)\in [0,1]\times [0,1]: (q-q^*(p))(p-p^*(q))<0\}.
\endaligned
\]
\end{theorem}
\begin{proof}
By Lemma \ref{lem:bistability-1}, {one has
\[
\Sigma_1=\{(p,q)\in [0,1]\times [0,1]: \Lambda(p,q)>0\}.
\]}
It then follows from Theorems {\ref{thm:Fpositive-1}, \ref{thm:Fpositive-1-a} and \ref{thm:Fpositive-2}} that $\Sigma_1$ is determined by  $(q-q^*(p))(q-p)<0$. The proof for {the characterization of} $\Sigma_2$ is {similar} and thus {skipped}. From the new characterizations for $\Sigma_1$ and $\Sigma_2$, it follows {that} $(p, q)\in \Sigma_3$ if and only if $(q-q^*)(q-p)>0$ and $(p-p^*)(p-q)>0$, which {amounts} to $(q-q^*(p))(p-p^*(q))<0$.
\end{proof}



\subsection{Stability of positive steady states of system \eqref{eq:main}}
The following result shows that any positive steady state of system \eqref{eq:main} is asymptotically stable. It is essentially due to He and Ni \cite{HeNi3}. For the sake of completeness, we have included here a slightly different demonstration of this result.

\begin{lemma}\label{bistability-3}
Suppose that {the free growth rate function} $a$ is non-constant and that $p$ {is not equal to} $q$. Then, any positive steady state of system \eqref{eq:main} is linearly stable and thus locally asymptotically stable.
\end{lemma}
\begin{proof} Let $(U,V)$ denote any positive steady state of {system} \eqref{eq:main}, \textit{i.e.}
\[ 
\left\{
\aligned
&D(p) U_{xx}+D(1-p) U_{yy}+(a-U-V)U=0 \quad \text{ in }\Omega,\\
&D(q) V_{xx}+D(1-q) V_{yy}+(a-U-V)V=0\quad\ \text{ in }\Omega,\\
&(D(p) U_x, D(1-p) U_y)\cdot \nu =(D(q) V_x, D(1-q) V_y)\cdot \nu=0\text{ on }\partial\Omega.
\endaligned
\right.
\] 

The linear stability of {this state} is determined by the sign of the principal eigenvalue $\lambda_1$ of the linear problem
\[ 
\left\{
\aligned
&D(p)\varphi_{xx}+D(1-p)\varphi_{yy}+(a-2U-V)\varphi-\psi V+\lambda_1\varphi=0\text{ in }\Omega,\\
&D(q)\psi_{xx}+D(1-q) \psi_{yy}-U\varphi+(a-U-2V)\psi+\lambda_1\psi=0\text{ in }\Omega,\\
& (D(p) \varphi_x, D(1-p) \varphi_y)\cdot \nu=0\text{ on }\partial\Omega,\\
&(D(q) \psi_x, D(1-q) \psi_y)\cdot \nu=0\text{ on }\partial\Omega.
\endaligned
\right.
\] 

It is known (see \cite{CC2003,Smith}) that we may choose $\varphi>0$ and $\psi<0$ in $\overline\Omega$. Set $W=\varphi/U$ and $Z=-\psi/V$ {so that $W$ and $Z$ are both positive} in $\overline\Omega$ and satisfy
\[ 
\left\{
\aligned
&D(p) (U^2 W_x)_x+D(1-p) (U^2 W_y)_y-U^3 W+U^2VZ+\lambda_1 U^2W=0\text{ in }\Omega,\\
&D(q) (V^2 Z_x)_x+D(1-q) (V^2 Z_y)_y + UV^2 W-V^3 Z+ \lambda_1 V^2 Z=0\text{ in }\Omega,\\
& (D(p) W_x, D(1-p) W_y)\cdot \nu=(D(q) Z_x, D(1-q) Z_y)\cdot \nu=0\text{ on }\partial\Omega.
\endaligned
\right.
\] 
Multiplying the first line {of this system} by $W^2$ and integrating the result {over $\Omega$ yields}
\[
\int_\Omega \left[D(p)U^2 W (W_x)^2+D(1-p)U^2 W (W_y)^2
+(UW)^3
-(UW)^2 (VZ)-\lambda_1 U^2 W^3\right]=0.
\]
Similarly, multiplying the {second line of the system} by $Z^2$ and
integrating the result {over $\Omega$}, we {find that}
\[
\int_\Omega \left[D(q)V^2 Z_x(Z^2)_x+D(1-q)V^2 Z_y (Z^2)_y+(VZ)^3-(UW)(VZ)^2-\lambda_1 V^2 Z^3\right]=0.
\]

It suffices to show $\lambda_1>0$. We argue by contradiction {by assuming} that $\lambda_1\le0$. Then, {one has}
\begin{equation}\label{eq:nihe-1}
\int_\Omega (UW)^3\le \int_\Omega (UW)^2 (VZ),
\end{equation}
and the equality in \eqref{eq:nihe-1} holds if and only if $\lambda_1=0$, $W$ is a positive constant, and $UW$ is a positive scalar multiple of $VZ$. Similarly, if $\lambda_1\le 0$, one has
\begin{equation}\label{eq:nihe-2}
\int_\Omega (VZ)^3\le \int_\Omega (UW)(VZ)^2,
\end{equation}
and the equality in \eqref{eq:nihe-2} holds if and only if $\lambda_1=0$, $Z$ is a positive constant, and $UW$ is a positive scalar multiple of $VZ$. {Finally, it follows from} the H\"older inequality that
\begin{equation}\label{eq:nihe-3}
\int_\Omega (UW)^3\le \left(\int_\Omega (UW)^3\right)^{2/3}\left( \int_\Omega (VZ)^3\right)^{1/3}
\end{equation}
and
\begin{equation}\label{eq:nihe-4}
\int_\Omega(VZ)^3
\le \left(\int_\Omega (VZ)^3\right)^{2/3}\left( \int_\Omega (UW)^3\right)^{1/3},
\end{equation}
from which we see that inequalities \eqref{eq:nihe-1}, \eqref{eq:nihe-2}, \eqref{eq:nihe-3} and \eqref{eq:nihe-4} must all be equalities. {As a consequence}, $\lambda_1$ {is zero}, both $W$ and {$Z$} are positive constants, and $UW$ is a positive scalar multiple of $VZ$, i.e. $U=cV$ for some positive constant {$c$}. Therefore, $U$ satisfies
\[ 
\left\{
\aligned
&D(p) U_{xx}+D(1-p) U_{yy}+(a-(c+1)U)U=0\quad \text{ in }\Omega,\\
&D(q) U_{xx}+D(1-q) U_{yy}+(a-(c+1)U)U=0\quad \text{ in }\Omega,\\
& (D(p) U_x,D(1-p) U_y)\cdot\nu=(D(q)U_x,D(1-q)U_y)\cdot\nu=0\text{ on }\partial\Omega.
\endaligned
\right.
\] 
Hence, {we find that} $U\equiv N_p/(c+1)$ and $U\equiv N_q/(c+1)$, which implies that $N_p\equiv N_q$. As {the function} $a$ is non-constant, {so is} $N_p$. Since $p\not= q$, by subtracting the equations {of the systems respectively satisfied by} $N_p$ and $N_q$, we see that $N_p$ is {a solution to system} \eqref{eq:wave-1} {and is therefore a non-constant function, which contradicts the assumption.}
\end{proof}

\subsection{Global dynamics of system \eqref{eq:main} 
}
As the two-species competition model \eqref{eq:main} is strongly monotone, its global dynamics can be fully determined by the local stability of {its} equilibria in some cases, see \cite[Chapter IV]{Hess} {for instance}. {Let us} recall {below} some {known facts}.

\begin{enumerate}[(a)]
\item If there is no positive steady state, then one of the semi-trivial equilibria is unstable and the other is globally asymptotically stable among non-negative and non-identically zero initial data.\label{global dynamics a}

\item If there is a unique positive steady state and it is stable, then it is globally asymptotically stable.\label{global dynamics b}

\item If all positive steady states are asymptotically stable, then there is at most one of them. In particular, either \eqref{global dynamics a} or \eqref{global dynamics b} applies.\label{global dynamics c}
\end{enumerate}

\medskip

We are now ready to infer on the global stability of steady states.

\begin{theorem}\label{thm:Fpositive-1-aa}
Suppose {that the function $F$ is positive} in $[0,1]$ {and} let $p^*$ and $q^*$ be the functions {introduced in} Theorems {\rm\ref{thm:Fpositive-1}} and {\rm\ref{thm:Fpositive-1-a}}, {respectively}. Then, {one of the following statements} holds.

\begin{enumerate}[\rm(i)]
\item If 
$p<q<q^*$, then {the steady state} $(N_p, 0)$ is globally asymptotically stable;\label{thm:Fpositive-1-aa - statement 1}

\item If $q<p<p^*$,
then {the steady state} $(0, N_q)$ is globally asymptotically stable;\label{thm:Fpositive-1-aa - statement 2}

\item If 
either $q^*<q\le 1$ or $p^*<p\le 1$ holds, then system \eqref{eq:main} has a unique positive steady state, which is also globally asymptotically stable among non-negative and not identically zero initial data.\label{thm:Fpositive-1-aa - statement 3}
\end{enumerate}
\end{theorem}
\begin{proof} We first establish statement \eqref{thm:Fpositive-1-aa - statement 1}. By Theorem \ref{thm:bistability-2}, we see that under assumption $p\in [0, 1]$ and $p<q<q^*$, $\Lambda(p, q)>0>\Lambda(q, p)$. Hence, $(N_p, 0)$ is stable and $(0, N_q)$ is unstable. As system \eqref{eq:main} is strongly monotone, by Lemma \ref{bistability-3} and statements \eqref{global dynamics a} and \eqref{global dynamics c}, $(N_p, 0)$ is globally stable.

The proof of {statement} \eqref{thm:Fpositive-1-aa - statement 2} is similar {to} that of {statement} \eqref{thm:Fpositive-1-aa - statement 1} and thus omitted.

For {statement} \eqref{thm:Fpositive-1-aa - statement 3}, $\Lambda(p, q)<0$ and $\Lambda(q, p)<0$. Hence, both states $(N_p,0)$  and $(0,N_q)$ are unstable. As system \eqref{eq:main} is strongly monotone, by Lemma \ref{bistability-3} and statements \eqref{global dynamics b} and \eqref{global dynamics c} recalled above, there is a unique positive steady state which is globally asymptotically stable.
\end{proof}



Note that $p^*\equiv 1$ if and only if $q^*\equiv 1$. For such scenario, alternative \eqref{thm:Fpositive-1-aa - statement 3} in Theorem \ref{thm:Fpositive-1-aa} does not occur, {the state} $(N_p, 0)$ is globally stable when $q>p$ and {the state} $(0, N_q)$ is globally stable when $q<p$.

Similar to Theorem \ref{thm:Fpositive-1-aa}, if $F<0$ in $[0, 1]$, the global dynamics of system \eqref{eq:main} can be characterized as follows.

\begin{theorem}\label{thm:Fpositive-2-aa}
Suppose {that the function $F$ is negative} in $[0,1]$ {and} let $p^*$ and $q^*$ be the functions {introduced in} Theorems {\rm\ref{thm:Fpositive-2}} and {\rm\ref{thm:Fpositive-2-a}}, {respectively}. Then, {one of the following statements} holds.
\begin{enumerate}[\rm(i)]
\item If 
$q^*<q<p$, then {the steady state} $(N_p,0)$ is globally asymptotically stable.

\item If $p^*<p<q$,
then {the steady state} $(0,N_q)$ is globally asymptotically stable.

\item If 
either $0\le q<q^*(p)$ or $0\le p<p^*(q)$ holds, then system \eqref{eq:main} has a unique positive steady state, which is also globally asymptotically stable among non-negative and not identically zero initial data.
\end{enumerate}
\end{theorem}

If {the function} $F$ changes sign in $(0, 1)$, the global dynamics of system \eqref{eq:main} is given by the following result:

\begin{theorem}\label{thm:Fchangessign-1}
Suppose {that the function} $F$ changes sign in $(0,1)$ {and} let $p=p^*(q)$ and $q=q^*(p)$ be {the} functions given in Theorems {\rm\ref{thm:Fpositive-3}} and {\rm\ref{thm:Fpositive-4}}. Then, {one of the following alternatives} holds.
\begin{enumerate}[\rm(i)]
\item If $(q-q^*)(p-p^*)<0$, then {system} \eqref{eq:main} has a unique positive steady state, which is also globally asymptotically stable among non-negative and not identically zero initial data.

\smallskip

\item If $(q-q^*)(q-p)<0$, then {the steady state} $(N_p,0)$ is globally asymptotically stable.

\smallskip

\item If $(q-q^*)(q-p)<0$, then {the steady state} $(0,N_q)$ is globally asymptotically stable.
\end{enumerate}
\end{theorem}

The proof of Theorem \ref{thm:Fchangessign-1} is the same as that of Theorem \ref{thm:Fpositive-1-aa} {and} follows from Theorem \ref{thm:bistability-2} and Lemma \ref{bistability-3}.

\section{Discussions}\label{sec:discussions}
In this paper, we considered a reaction-diffusion model for two competing populations, which disperse in a bounded two  dimensional habitat {by moving} horizontally and vertically with different {probabilities} but  are otherwise identical. We regard these probabilities as dispersal strategies and ask what strategies are evolutionarily stable.
 
Our main finding is  that the only evolutionarily stable dispersal strategies are to move in one direction. In particular, when the resources are distributed inhomogeneously only in one direction, \textit{e.g.} horizontally, our result implies that the evolutionarily stable strategy could simply be to move in the vertical direction, in which the resources are homogeneously distributed. More precisely, we introduced a function $F$ of the dispersal probability, which measures the difference between the spatial variations of the population equilibrium  distributions in horizontal and vertical directions: when it is positive, the species has more  variations in the horizontal direction; when it is negative, it has more variations in the vertical direction. We show that function $F$ is monotone decreasing and  that the  evolutionarily stable dispersal strategies are to maximize the function $F$ when it is positive and to minimize it when it is negative. 
As the population distribution at equilibrium is often positively  correlated with the resource distribution, function $F$ also indirectly measures the difference between the resource variations in horizontal and vertical directions. Therefore, our results seem to predict that it is more favorable for the species to choose the direction with smaller variations in resource distributions. This finding seems to be in agreement with the classical results of Hasting \cite{Hastings1983} and Dockery \textit{et al.} \cite{DHMP} for the evolution of slow dispersal, \textit{i.e.}, random diffusion is selected against in spatially heterogeneous  environments.

We further investigated the local and global dynamics of the two-species system and determined the  dynamics of system \eqref{eq:main} for three different cases of the selection gradient. 
We applied numerical simulations to illustrate how the shapes of function $F$, the local stability of the semi-trivial steady states and the global dynamics of the system sequentially change with respect to {a} certain parameter which measures the difference between the resource variations in {the} horizontal and vertical directions. Our numerical results  suggest that if the spatial variations of resource distributions in  vertical and horizontal directions become more comparable, the chances for the coexistence of two competing populations could be greater.

One of our future works is to extend the mathematical modelling and analysis to any dimensional habitats, and to continuous trait models. Another future work will be to include a temporal variation of the environment and ask how it affects the evolution of horizontal and vertical movement. For example, if we choose $a(t,x,y)=\lambda A(x)+(1-\lambda)B(t,y)$, a natural question is when vertical movement will be selected as in \cite{HMP}.

\bigskip

\noindent{\bf Acknowledgement.}
We thank Professor Yoshikazu Giga for the helpful discussions which motivated the study of anisotropic diffusion. We are also very grateful to Maxime Chupin for his help with the post-processing of our numerical simulations.
\bigskip

\bibliographystyle{plain}
\bibliography{blls}

\section{Appendix}

\subsection{Some remarks on {solutions to a} wave equation}
In the proofs of Lemmas \ref{lem:selection-gradient-1} and \ref{lemma:concavity},
the following result, which seems to be of self interest,  plays an important role in eliminating the degeneracy of the function $F$ and in establishing the strict concavity of the function $\Lambda(p,q)$ with respect to $q$: 
\begin{lemma}\label{lemma:wave}
Let $W$ {in} $C^2(\Omega)\cap C^1(\overline\Omega)$ be a solution to {the system}
\begin{equation}\label{eq:wave-1}
\left\{
\aligned
&W_{xx}-W_{yy}=0\text{ in }\Omega,\\
&W_x\nu_x=W_y\nu_y=0\text{ on }\partial\Omega.
\endaligned
\right.
\end{equation}
Then, {the function} $W$ {is} constant.
\end{lemma}
\begin{proof}
{By the strict convexity} assumption {on the domain, the components} $\nu_x$ and $\nu_y$ {of the outward normal vector $\nu$} are non-zero on {the boundary} $\partial\Omega$, except {possibly over} a set of measure zero. Hence, $W_x$ and $W_y$ {both vanish almost everywhere} on $\partial\Omega$. Since $W$ {belongs to} $C^1(\overline\Omega)$, {the gradient} $\nabla W$ {vanishes} on $\partial\Omega$.

Set $\eta=x+y$, $\zeta=x-y$ and $Z(\eta,\zeta):=W(x,y)$. {The function $Z$ then} satisfies
\[
Z_{\eta\zeta}=0\text{ in }\Omega'\quad \text{ and }\quad {(Z_\eta,Z_\zeta)=(0,0)}\text{ on }\partial\Omega',
\]
where $\Omega'$ is the image of $\Omega$ {under the map} $(x,y)\mapsto(\eta,\zeta)$. {It follows from the first relation that}  $Z(\eta,\zeta)=f(\eta)+g(\zeta)$ for some functions $f$ and $g$, {and the second one then} implies that both $f$ and $g$ {have to} be constant functions. {As a consequence,} $Z$ is a constant function, {and so} is $W$.
\end{proof}


It is possible to construct domains such that {problem} \eqref{eq:wave-1} admits \emph{non-constant} solutions, if we allow Lipschitz domains with flat parts on their boundaries.

\medskip
\noindent{\bf Example 1.}
Consider $\Omega=(0, 1)\times (0, 1)$ {and} let $f$ be an even
and $2$-periodic function in $\mathbb{R}$. Set $W(x,y)=f(x+y)+f(x-y)$, {which then clearly} satisfies {problem} \eqref{eq:wave-1}, {and is a positive non-constant function} if $f$ is taken positive and non-constant.

\medskip

On the other hand, the type of {domain} given in {the above example seems} to be non-generic, as illustrated by the following result:
\begin{lemma}
Suppose that $\Omega=(0, L_1)\times(0, L_2)$ for {some} positive numbers $L_1$ and $L_2$. If $L_1/L_2$ is not a rational number, then {problem} \eqref{eq:wave-1} has only constant solutions.
\end{lemma}

\begin{proof}
For any $W$ {satisfying problem} \eqref{eq:wave-1}, we have $W(x,y)=f(x+y)+f(x-y)$ for some scalar function $f$. Then $W_x=f'(x+y)+f'(x-y)$. {Since} $W_x(0, y)=0$, we have $f'(y)=-f'(-y)$, i.e. $f'$ is an odd function.
{Since} $W_x(L_1, y)=0$, we have $f'(y+L_1)=-f'(L_1-y)=f'(y-L_1)$, i.e. $f'$ is $2L_1$-periodic.

Similarly, $W_y=f'(x+y)-f'(x-y)$. Note that $W_y(x,0)=0$ automatically holds. By $W_y(x,L_2)=0$, we have $f'(x+L_2)=f'(x-L_2)$, i.e. $f'$ is also $2L_2$-periodic. Hence, if $L_1/L_2$ is not rational, then $f'$ must be a constant function. Since $f'$ is an odd function, then $f'=0$, i.e. $W$ is a constant function.
\end{proof}

\subsection{A remark about {a} possible {degeneracy induced by the domain} $\Omega$}\label{appendix:2} 
Throughout the paper, we have assumed that $\Omega$ is a strictly convex domain. We now comment on this {point}, showing with a very basic example that a domain with flat {parts} on {its} boundary may lead to {a} degeneracy of {the} function $F$.

Consider $\Omega=(0, 1)\times (0, 1)$.
Let $f$ and $D$ be given in Example 1.
Set
$$
a:=-(\overline{D}+\underline{D})\frac{W_{xx}}{W}+W.
$$
It is easy to check that
for each $\theta\in [0, 1]$,
$W$ also solves \eqref{eq:singlespeciessteady},
i.e. $N_\theta\equiv W$ for each $\theta\in [0, 1]$.
Since $W$ is non-constant, $a$ is also non-constant.
Furthermore,  for each $p, q\in [0, 1]$,
\eqref{eq:main} has a {continuum} of positive steady states
of the form $(U, V)=(sW, (1-s)W)$, $s\in (0, 1)$.
Moreover, in this case,
$F\equiv 0$ in $[0, 1]$ even though the function $a$ is not constant. Indeed,
since $f$ is even, we have $W(x,y)=W(y,x)$, which implies that $F\equiv 0$ in $[0, 1]$.

\subsection{The sign of $\Lambda(0,1)$}\label{appendix:3}
In this subsection we construct  an example 
to support the claim
made in Remark \ref{remark:11},
as shown  by the following result: 
 
\begin{proposition}\label{prop:1}
Assume that $a(x,y)=\lambda A(x)+(1-\lambda)A(y)$, where $\lambda$ {belongs to} $[0,1)$, and $A$ is positive H\"older continuous function, which attains a strict global maximum at $y=\hat{y}$ for some $\hat{y}$, \textit{i.e.} $A(y)<A(\overline{y})$ for every $y\not=\hat{y}$. Then, there exists some positive real number $\delta$ such that for  $\overline{D}>1/\delta$, then for sufficiently small $\underline{D}$, $\Lambda(0, 1)<0$.
\end{proposition}

We first establish some a {\it priori} estimate on $N_0$. 
By definition, $N_0$ satisfies
\[ 
\left\{
\aligned
&\underline{D}{(N_{0})_{xx}}+\overline{D}{(N_{0})_{yy}}+(a-N_0)N_0=0\text{ in }\Omega,\\
&(\underline{D}{(N_{0})_{x}},\overline{D}{(N_{0})_{y}})\cdot\nu=0\text{ on }\partial\Omega.
\endaligned
\right.
\] 

Define $\Omega_x:=\{y: (x, y)\in \Omega\}$ and $\Omega_y:=\{x:(x, y)\in \Omega\}$,  and let $\Omega:=\cup_{y_*<y<y^*} \Omega_y$.  If $\Omega$ is strictly convex, we may express it as
\[
\Omega=\left\{(x,y):  x_*<x< x^*,\ y^*(x)<y<y_*(x)\right\}
\]
for some $x_*<x^*$,  $y_*(x)<y^*(x)$ in $(x_*, x^*)$ and $y_*(x)=y^*(x)$ at $x=x_*, x^*$. For any $x\in (x_*, x^*)$, clearly $\Omega_x:=\{(x, y): y_*(x)<y<y^*(x)\}$.

\begin{lemma}\label{lemma:temp011}
Suppose that the function $a$ satisfies {\rm(A1)}, and $\Omega$ is strictly convex, $C^1$. For any $\epsilon>0$, there exists $\delta>0$ such that if  $\overline{D}>1/\delta$, then for sufficiently small $\underline{D}$,
\[
\frac{1}{|\Omega_x|} \int_{\Omega_x} a(x, y)\,dy-\epsilon
\le N_0(x, y)\le \frac{1}{|\Omega_x|} \int_{\Omega_x} a(x, y)\,dy+\epsilon
\]
holds for any $(x, y)\in\Omega$.
\end{lemma}
\begin{proof} Given $\epsilon>0$ small, choose function $\tilde{a}\in C^2(\overline\Omega)$ such that $a+\epsilon<\tilde{a}\le a+2\epsilon$ in $\Omega$ and $\tilde{a}(x,y)$ is constant 
for $x\in [x_*, x_*+\epsilon]\cup [x^*-\epsilon, x^*]$ and  $y_*(x)<y<y^*(x)$. For each $x\in (x_*, x^*)$, let $\tilde{N}(x, y)$ denote the unique positive solution of the equation
\[
\left\{
\aligned
&\overline{D}\tilde{N}_{yy}+\tilde{N}(\tilde{a}-\tilde{N})=0\quad \text{ in }y_*(x)<y<y^*(x),\\
&\tilde{N}_y(y_*(x))=-\frac{\epsilon}{\overline{D}},\quad \tilde{N}_y(y^*(x))=\frac{\epsilon}{\overline{D}}.
\endaligned
\right.
\]
Then $\tilde{N}$ is independent of $x$ for $x\in [x_*, x_*+\epsilon]\cup [x^*-\epsilon, x^*]$, for which it satisfies $\tilde{N}_x=0$. Therefore, $\tilde{N}$ satisfies
\[
\underline{D}\tilde{N}_x\nu_x+\overline{D}\tilde{N}_y\nu_y
=\epsilon|\nu_y|
\ge 0
\]
for $(x, y)\in\partial\Omega$ and $x\in [x_*, x_*+\epsilon]\cup [x^*-\epsilon, x^*]$. For $(x, y)\in\partial\Omega$ and $x\in [x_*+\epsilon, x^*-\epsilon]$,
\[
\underline{D}\tilde{N}_x\nu_x+\overline{D}\tilde{N}_y\nu_y
\ge -\underline{D}\|\tilde{N}_x\nu_x\|_{L^\infty}+\epsilon \min_{x\in [x_*+\epsilon, x^*-\epsilon]}|\nu_y|>0
\]
for sufficiently small $\underline{D}$,
since $|\nu_y|$ is a strictly positive and continuous function for $(x, y)\in\partial\Omega$ and $x\in [x_*+\epsilon, x^*-\epsilon]$. This implies that, for sufficiently small $\underline{D}$, $\underline{D}\tilde{N}_x\nu_x+\overline{D}\tilde{N}_y\nu_y\ge 0$ holds on $\partial\Omega$.

In $\Omega$, $\tilde{N}$ satisfies
\[
\underline{D}\tilde{N}_{xx}+\overline{D}\tilde{N}_{yy}+\tilde{N}({a}-\tilde{N})
=\underline{D}\tilde{N}_{xx}+\tilde{N}({a}-\tilde{a})
\le \underline{D}\tilde{N}_{xx}-\tilde{N}\epsilon<0,
\]
provided that $\underline{D}$ is sufficiently small. This implies that $\tilde{N}$ is a super-solution for the equation of $N_0$. Hence, if $\underline{D}$ is small, then $N_0\le \tilde{N}$ in $\Omega$. As $\overline{D}$ tends to infinity, $\tilde{N}$ converges uniformly to $\frac{1}{|\Omega_x|} \int_{\Omega_x} \tilde{a}(x, y)\,dy$ in $\Omega$. Hence, there exists some $\delta>0$ such that for $\overline{D}>1/\delta$, if $\underline{D}$ is small, then
\[
\tilde{N}\le \frac{1}{|\Omega_x|} \int_{\Omega_x} \tilde{a}(x, y)\,dy+\epsilon,
\]
which implies that
\[
N_0\le \tilde{N}\le \frac{1}{|\Omega_x|} \int_{\Omega_x} {a}(x, y)\,dy+3\epsilon
\]
holds in $\Omega$. The lower bound of $N_0$ can be similarly established.
\end{proof}

Set $\Omega:=\cup_{\underline{y}<y<\overline{y}} \Omega_y$ for some $\underline{y}<\overline{y}$.

\begin{lemma}\label{lemma:eigenvalue-estimate-a}
The following result holds:
\begin{equation}\label{eq:eigenvalue-estimate-a}
\limsup_{\overline{D}\to \infty}\limsup_{\underline{D}\to 0}\Lambda(0,1)\le\min_{\underline{y}\le y\le \overline{y}}\frac{1}{|\Omega_y|}\int_{\Omega_y}
\left(\frac{1}{|\Omega_x|}\int_{\Omega_x} a(x,z)\,dz-a(x,y)\right)\,dx.
\end{equation}
\end{lemma}
In particular, there exists some $\delta>0$ such that for  $\overline{D}>1/\delta$, then for sufficiently small $\underline{D}$, $\Lambda(0, 1)<0$, provided that
\[
\int_{\Omega_y}\left(\frac{1}{|\Omega_x|}\int_{\Omega_x} a(x,z)\, dz-a(x,y)\right)\, dx\not\equiv 0\qquad \mbox{for}\ y\in [\underline{y}, \overline{y}].
\]
\begin{proof}
By Lemma \ref{lemma:temp011}, for any $\epsilon>0$, there exists some $\delta>0$ such that for  $\overline{D}>1/\delta$, then for sufficiently small $\underline{D}$, $N_0(x,y)\le \frac{1}{|\Omega_x|}\int_{\Omega_x} a(x,z)\, dz+\epsilon$ in $\Omega$.

Note that
\[
\Lambda(0,1)=\inf_{\varphi\not=0,\ \varphi\in H^1(\Omega)}\frac{\int_\Omega[\overline{D}(\varphi_x)^2+\underline{D}(\varphi_y)^2-(a-N_0)\varphi^2]}{\int_\Omega \varphi^2}.
\]
By choosing $\varphi=\varphi(y)$ we have
\[
\aligned
\Lambda(0,1)
&\le\frac{\int_\Omega [\underline{D}(\varphi_y)^2-(a-N_0)\varphi^2]}{\int_\Omega\varphi^2}\\
&\le\frac{\int_\Omega [\underline{D}(\varphi_y)^2-(a-\frac{1}{|\Omega_x|}\int_{\Omega_x} a(x,z)\, dz-\epsilon)\varphi^2]}{\int_\Omega \varphi^2},
\endaligned
\]
which implies that
\[
\aligned
\limsup_{\underline{D}\to 0} \Lambda(0,1)
&\le-\frac{\int_\Omega(a-\frac{1}{|\Omega_x|}\int_{\Omega_x} a(x,z)\, dz)\varphi^2}{\int_\Omega \varphi^2}+\epsilon\\
&:=\frac{\int_{\underline{y}}^{\overline{y}} G(y)\varphi^2(y)|\Omega_y|\, dy}{ \int_{\underline{y}}^{\overline{y}}\varphi^2(y)|\Omega_y|\,dy}
+\epsilon
\endaligned
\]
holds for any $\epsilon>0$ and $\varphi(y)\not=0$, where
\[
G(y):=\frac{1}{|\Omega_y|}\int_{\Omega_y}
\left(\frac{1}{|\Omega_x|}\int_{\Omega_x} a(x,z)\, dz-a(x,y)\right)\, dx.
\]
Hence, given any $x$ and taking the infimum of the right-hand side of the above inequality over all $\varphi=\varphi(y)\not=0$, we have
\[
\limsup_{\overline{D}\to\infty}\limsup_{\underline{D}\to 0}\Lambda(0,1)
\le \min_{\underline{y}\le y\le \overline{y}} G(y)+\epsilon.
\]
By passing $\epsilon\to 0$ we see that \eqref{eq:eigenvalue-estimate-a} holds.

Finally, we claim that $\int_{\underline{y}}^{\overline{y}} G(y)|\Omega_y|\, dy=0$:
\[
\aligned
\int_{\underline{y}}^{\overline{y}} G(y)|\Omega_y|\, dy
&=\int_{\underline{y}}^{\overline{y}}\int_{\Omega_y}\int_{\Omega_x} \frac{a(x,z)}{|\Omega_x|}\, dz\, dx\,dy-\int_\Omega a(x,y)\,dx\,dy\\
&=\int_{\underline{y}}^{\overline{y}}\int_{\Omega_z}\int_{\Omega_x} \frac{a(x,z)}{|\Omega_x|}\, dy\, dx\,dz-\int_\Omega a(x,y)\,dx\,dy
\\
&=\int_{\underline{y}}^{\overline{y}}\int_{\Omega_z}a(x,z)\, dx\,dz-\int_\Omega a(x,y)\,dx\,dy=0.
\endaligned
\]
Hence, if $G(y)\not\equiv$ 0, we have $\min_{\underline{y}\le y\le \overline{y}} G(y)<0$. This completes the proof.
\end{proof}

\begin{proof}[Proof of Proposition {\rm\ref{prop:1}}]
If $a(x, y)=\lambda A(x)+(1-\lambda)A(y)$, then
\[
G(y)=(1-\lambda)\left[\int_{\Omega_y} \frac{\int_{\Omega_x}A(z)\, dz}{|\Omega_x|}\,dx-A(y)|\Omega_y|\right].
\]
To apply Lemma \ref{lemma:eigenvalue-estimate-a}, it suffices to check
\[
\frac{1}{|\Omega_y|}\int_{\Omega_y} \frac{\int_{\Omega_x}A(z)\, dz}{|\Omega_x|}\,dx\not\equiv A(y).
\]
To see this, observe that for all $y\not=\hat{y}$ and $y\in (\underline{y}, \overline{y})$,
\[
\frac{1}{|\Omega_y|}\int_{\Omega_y} \frac{\int_{\Omega_x}A(z)\, dz}{|\Omega_x|}\,dx<\frac{1}{|\Omega_y|}\int_{\Omega_y} \frac{\int_{\Omega_x}A(\hat{y})\, dz}{|\Omega_x|}\,dx
=A(\hat{y}).
\]
The conclusion thus follows from Lemma \ref{lemma:eigenvalue-estimate-a}.
\end{proof}

\end{document}